# Optimal control of Goursat-Darboux systems in domains with curvilinear boundaries

S. A. Belbas
Mathematics Department
University of Alabama
Tuscaloosa, AL. 35487-0350. USA.

e-mail: SBELBAS@GP.AS.UA.EDU

Abstract. We derive necessary conditions for optimality in control problems governed by hyperbolic partial differential equations in Goursat-Darboux form. The conditions consist of a set of Hamiltonian equations in Goursat form, side conditions for the Hamiltonian equations, and an extremum principle akin to Pontryagin's maximum principle. The novel ingredient is that the domain over which the optimal control problem is posed is not rectangular. The non-rectangular nature of the domain affects the optimality conditions in a substantial way.

## 1. Introduction.

The theory of optimal control for Goursat-Darboux systems is an important extension of the corresponding theory for controlled ordinary differential equations. This extension has the same place in mathematical optimal control theory that the extension from single integrals to multiple integrals has in the calculus of variations.
Generally, a controlled Goursat equation has the form

$$\frac{\partial^2 x(s,t)}{\partial s \partial t} = f(s,t,x(s,t),x_s(s,t),x_t(s,t),u(s,t)) \qquad \text{--- (1.1)}$$

Additionally, we need to specify a domain, say G, for the two-dimensional variable (s, t), and appropriate boundary conditions on $\partial G$. The simplest form of a domain G is a rectangle, $G=\{(s,t): 0 \le s \le a, 0 \le t \le b\}$. In that case, we have the standard Goursat-Darboux problem with data on characteristics, i.e. (1.1) is accompanied by the set of boundary conditions

$$\begin{aligned} &x(s,0) = x_1(s), \text{ for } 0 \le s \le a; \\ &x(0,t) = x_2(t), \text{ for } 0 \le t \le b; \\ &x(0,0) = x_0 \equiv x_1(0) = x_2(0) \end{aligned} \qquad \text{--- (1.2)}$$

It is well known (e.g. [CC]) that, under natural Lipschitz conditions on the function f, continuity of the boundary data $x_1$ and $x_2$, and the consistency condition $x_1(0) = x_2(0)$, the problem {(1.1), (1.2)} has a unique solution. There are also more general existence and uniqueness results, but those are not relevant to the optimal control issues that we study in the present paper.
An optimal control problem for the system {(1.1), (1.2)} concerns the minimization of a functional J given by

$$\begin{aligned} J := &\iint_G F(s,t,x(s,t),x_s(s,t),x_t(s,t),u(s,t))dA(s,t) + \\ &+ \int_0^a F_1(s,x(s,b),x_s(s,b))ds + \int_0^b F_2(t,x(a,t),x_t(a,t))dt + F_0(x(a,b)) \end{aligned} \qquad \text{--- (1.3)}$$

Necessary conditions for optimality, in a form analogous to Pontryagin's maximum principle, have be obtained for the problem {(1.1), (1.2), (1.3)} in [BDMO, E1, E2, E3,

PS, S, VST]; related work, from the point of view of dynamic programming with two-dimensional "time" variable, has been done in [B1, B2, B3].

We define the Hamiltonian

$$H(s,t,x,p,q,\psi,u) := F(s,t,x,p,q,\psi,u) + \psi f(s,t,x,p,q,u) \qquad \text{--- (1.4)}$$

(The variables p, q stand for the slots of $x_s$, $x_t$, respectively.) Under conditions of sufficient differentiability of all functions involved, the Hamiltonian equations for the co-state $\psi(s, t)$ are

$$\frac{\partial^2\psi(s,t)}{\partial s\partial t} = \frac{\partial H(s,t,x(s,t),x_s(s,t),x_t(s,t),\psi(s,t),u(s,t))}{\partial x} -$$
$$-\frac{D}{Ds}\left(\frac{\partial H(s,t,...)}{\partial p}\right) - \frac{D}{Dt}\left(\frac{\partial H(s,t,...)}{\partial q}\right), \text{ for } (s,t) \in \text{int}(G);$$
$$\frac{\partial\psi(s,b)}{\partial s} + \frac{\partial H(s,b,...)}{\partial q} - \frac{D}{Ds}\left(\frac{\partial F_1(s,...)}{\partial p}\right) + \frac{\partial F_1(s,...)}{\partial x} = 0, \text{ for } 0 \le s < a;$$
$$\frac{\partial\psi(a,t)}{\partial t} + \frac{\partial H(a,t,...)}{\partial p} - \frac{D}{Dt}\left(\frac{\partial F_2(t,...)}{\partial q}\right) + \frac{\partial F_2(t,...)}{\partial x} = 0, \text{ for } 0 \le t < b;$$
$$-\psi(a,b) + \frac{\partial F_1(a,...)}{\partial p} + \frac{\partial F_2(b,...)}{\partial q} + \frac{\partial F_0(y(a,b))}{\partial y} = 0 \qquad \text{--- (1.5)}$$

The operators $\frac{D}{Ds}, \frac{D}{Dt}$ are the operators of taking total derivatives with respect to s, t, respectively. This means that for any function $\Phi(s,t,x(s,t),x_s(s,t),x_t(s,t),\psi(s,t),u(s,t))$, we have

$$\frac{D\Phi}{Ds} = \frac{\partial\Phi}{\partial s} + \frac{\partial\Phi}{\partial x}\frac{\partial x}{\partial s} + \frac{\partial\Phi}{\partial p}\frac{\partial^2 x}{\partial s^2} + \frac{\partial\Phi}{\partial q}\frac{\partial^2 x}{\partial s\partial t} + \frac{\partial\Phi}{\partial\psi}\frac{\partial\psi}{\partial s} + \frac{\partial\Phi}{\partial u}\frac{\partial u}{\partial s};$$
$$\frac{D\Phi}{Dt} = \frac{\partial\Phi}{\partial t} + \frac{\partial\Phi}{\partial x}\frac{\partial x}{\partial t} + \frac{\partial\Phi}{\partial p}\frac{\partial^2 x}{\partial s\partial t} + \frac{\partial\Phi}{\partial q}\frac{\partial^2 x}{\partial t^2} + \frac{\partial\Phi}{\partial\psi}\frac{\partial\psi}{\partial t} + \frac{\partial\Phi}{\partial u}\frac{\partial u}{\partial t} \qquad \text{--- (1.6)}$$

It is worth noticing that the terms $\frac{D}{Ds}\left(\frac{\partial H}{\partial p}\right)$ and $\frac{D}{Dt}\left(\frac{\partial H}{\partial q}\right)$ are peculiar to Goursat problems, and there are no analogous terms in the Hamiltonian equations for optimal control problems for systems governed by ordinary differential equations.

The condition of the existence of the derivatives $\frac{D}{Ds}\left(\frac{\partial H}{\partial p}\right)$ and $\frac{D}{Dt}\left(\frac{\partial H}{\partial q}\right)$ is a strong assumption, since, according to (1.6), it requires, among other things, differentiability of the control u with respect to s and t. One variant of the standard model of controlled Goursat-Darboux systems can be obtained by taking a control function u in the form

$$u(s,t) = \int_0^s \int_0^t v(\sigma,\tau)\, d\tau\, d\sigma \qquad \text{--- (1.7)}$$

so that, in reality, v becomes the control function; this approach has been introduced in [B1].

In cases in which it is not assumed that the functions appearing in the Hamiltonian equations for Goursat systems are sufficiently many times differentiable, it is possible to use an integral formulation of the Hamiltonian equations. The integral equation for the co-state then becomes

$$\psi(s,t) = \psi(a,t) + \psi(s,b) - \psi(a,b) + \int_s^a \int_t^b \frac{\partial H(\sigma,\tau,...)}{\partial x} d\tau\, d\sigma +$$
$$+ \int_s^a \left[\frac{\partial H(\sigma,t,...)}{\partial q} - \frac{\partial H(\sigma,b,...)}{\partial q}\right] d\sigma + \int_t^b \left[\frac{\partial H(s,\tau,...)}{\partial p} - \frac{\partial H(a,\tau,...)}{\partial p}\right] d\tau \ ;$$
$$\psi(s,b) - \psi(a,b) + \int_s^a \left[\frac{\partial H(\sigma,b,...)}{\partial q} + \frac{\partial F_1(\sigma,...)}{\partial x}\right] d\sigma + \frac{\partial F_1(s,...)}{\partial p} - \frac{\partial F_1(a,...)}{\partial p} = 0 \ ;$$
$$\psi(a,t) - \psi(a,b) + \int_t^b \left[\frac{\partial H(a,\tau,...)}{\partial p} + \frac{\partial F_2(\tau,...)}{\partial x}\right] d\tau + \frac{\partial F_2(s,...)}{\partial q} - \frac{\partial F_2(a,...)}{\partial q} = 0 \ ;$$
$$-\psi(a,b) + \frac{\partial F_1(a,...)}{\partial p} + \frac{\partial F_2(b,...)}{\partial q} + \frac{\partial F_0(x(a,b))}{\partial x} = 0 \qquad \text{--- (1.8)}$$

We have mentioned these two modifications, i.e. state dynamics with nonlocal operators acting on the control and integral formulation of the Hamiltonian equations, for the sake of completeness. We shall not use those modifications for the problems we analyze in the present paper, since the main feature of our work is the effect of the non-rectangular boundary on the optimality conditions. It is possible, however, to include the two modifications mentioned above in the problems over non-rectangular domains.

There are many reasons for studying optimal control problems for Goursat-Darboux equations over non-rectangular domains:

(1). A hyperbolic control problem may have been originally formulated in non-characteristic coordinates, over a domain that is rectangular in the original coordinates of the problem. When we change to characteristic coordinates, in order to take advantage of the techniques that have been developed for the optimal control of Goursat-Darboux systems, the domain will not generally be rectangular in the new coordinates.

(2). The techniques we develop, in this paper, for Goursat-Darboux problems over non-rectangular domains, can also be used for problems that involve certain integral terms, in the cost functional, with integrals supported on interior curves of the domain, as well as terms concentrated on interior points of the domain.

(3). There are also other reasons, in addition to relevance to applications, for studying Goursat-Darboux control problems over non-rectangular domains; these are "internal" reasons, they are intrinsic to the discipline of Mathematics. The study of Goursat-Darboux control problems over rectangular domains is not a mathematically satisfactory extension of the corresponding control theory for ordinary differential equations, because it does not make full use of the two-dimensional parameter $(s, t)$, which plays the role of generalized "time" for Goursat-Darboux systems. A crucial feature of two-dimensional space, compared with the one-dimensional line of real numbers, is the greater variety of the shapes of connected open sets that are possible on the plane. This feature is not exploited when attention is restricted to rectangular domains. Certain aspects of the mathematical structure that emerges from the intrinsic nature of variational calculus for controlled systems governed by Goursat-Darboux equations are obscured when attention is restricted to rectangular domains, and are revealed only when we examine more general non-rectangular domains.

2. Statement of the problem.

We consider a domain G in the first quadrant ($s \ge 0$, $t \ge 0$) of the st-plane, bounded by the segments [0, a] on the s-axis, the segment [0, b] on the t-axis, and a continuous and piecewise $C^2$ curve $(\gamma)$ with at most a finite number of points at which the normal vector field to $(\gamma)$ has jump discontinuities. The curve $(\gamma)$ is assumed to be a non-increasing graph, and every smooth segment of $(\gamma)$ can be represented in at least one of two forms: $s = \vartheta_1(t)$, $t = \vartheta_2(s)$. Clearly, both representations are possible on every smooth part of $(\gamma)$ that has nowhere a tangent parallel to any of the coordinate axes. We denote by $\mu$ the parameter of arc length on $(\gamma)$, measured from the point $P_0 \equiv (a,0)$ on the s-axis. The curve $(\gamma)$ is also parametrized by arc length, as $s = L_1(\mu)$, $t = L_2(\mu)$, $0 \le \mu \le L$, where L is the total length of $(\gamma)$. We consider a set **P** of points on $(\gamma)$, $\mathbf{P} := \{(a,0) \equiv P_0, P_1, P_2, ..., P_N, P_{N+1} \equiv (0,b)\}$, that contains all the points of jump discontinuities of the normal vector field on $(\gamma)$ but may also contain additional points. The reason for allowing the possibility of including additional points in **P** will become clear farther down in this section. We shall call the points in **P** vertices. We set $(\gamma_R) := (\gamma) \setminus \mathbf{P}$, and we shall call $(\gamma_R)$ the regular part of $(\gamma)$. For later use, we also define the regular part of G, denoted by $G_R$, as the set obtained by removing from G all straight lines, parallel to either the s-axis or the t-axis, and passing through a vertex. If $s(P_r)$, $t(P_r)$ are the s- and t- coordinates of each vertex $P_r$, then

$$G_R := \{(s,t) \in G : s \ne s(P_r) \text{ and } t \ne t(P_r) \ \forall r = 0,1,...,N+1 .$$

The fact that the curve $(\gamma)$ is a non-increasing curve implies the well-posedness, under standard conditions of Lipschitz continuity of f relative to x, $x_s$, $x_t$, and continuity of f with respect to u, of the Goursat-Darboux problem

$$\frac{\partial^2 x(s,t)}{\partial s \partial t} = f(s,t,x(s,t),x_s(s,t),x_t(s,t),u(s,t)), \text{ for } (s,t) \in (\text{int}\, G) \cup (\gamma);$$
$$x(s,0) = x_1(s) \text{ for } 0 \le s \le a; \ x(0,t) = x_2(t) \text{ for } 0 \le t \le b \qquad \text{--- (2.1)}$$

where the functions $x_1(.)$, $x_2(.)$ are continuous on $[0,a]$, $[0,b]$, respectively, and satisfy the consistency condition $x_1(0) = x_2(0)$, and we set $x_0 := x_1(0) = x_2(0)$. This follows from the fact that, under the stated conditions, the sets W(s, t), defined by

$$W(s,t) := \{(\sigma,\tau) : 0 \le \sigma \le s, \ 0 \le \tau \le t\} \qquad \text{--- (2.2)}$$

satisfy $W(s,t) \subseteq (\operatorname{int} G) \cup (\gamma)$ $\forall (s,t) \in (\operatorname{int} G) \cup (\gamma)$ , and the Goursat-Darboux problem (2.1) can be written in integral form as

$$x(s,t) = x_1(s) + x_2(t) - x_0 + \iint\limits_{W(s,t)} f(\sigma, \tau, x(\sigma,\tau), x_\sigma(\sigma,\tau), x_\tau(\sigma,\tau), u(\sigma,\tau))$$

for $(s,t) \in (\operatorname{int} G) \cup (\gamma)$

--- (2.3)

and well-posedness can be shown by standard techniques for such two-dimensional Volterra equations. We note that more general results of existence and uniqueness, under conditions that do not require the Lipschitz continuity of f with respect to all the variables we listed previously, are known in the research literature, but those results are not relevant to the control problems we study in this paper.

Thus, we consider a controlled dynamical system of the form (2.1), in which the state $x(s,t)$ takes values in n-dimensional Euclidean space, and the control function $u(s,t)$ takes values in m- dimensional Euclidean space. Coordinates will be denoted by supesrscripts, and the convention of tensor algebra about summation with respect to repeated subscripts and superscripts will be strictly used, except in cases in which we explicitly state that no summation is taken.

Before we define the cost functional that is to be minimized, we need some notation. We denote by $\vec{n}(s,t) = \begin{bmatrix} n_1(s,t) \\ n_2(s,t) \end{bmatrix}$ the unit normal vector field on the smooth part of $(\gamma)$ oriented towards the exterior of G. For any differentiable function $\varphi$, we denote by $\varphi_\mu$ or $\frac{\partial \varphi}{\partial \mu}$ the tangential derivative of $\varphi$ on $(\gamma)$, i.e.

$\varphi_\mu(s,t) = -n_2(s,t)\frac{\partial \varphi(s,t)}{\partial s} + n_1(s,t)\frac{\partial \varphi(s,t)}{\partial t}$ , for $(s,t) \in (\gamma_R)$. The outward normal derivative of $\varphi$ will be denoted by $\varphi_n$ or $\frac{\partial \varphi}{\partial n}$ , thus

$$\varphi_n(s,t) = n_1(s,t)\frac{\partial \varphi(s,t)}{\partial s} + n_2(s,t)\frac{\partial \varphi(s,t)}{\partial t}.$$

For later use, we define a third derivative, denoted by $\varphi_{n^*}$ or $\frac{\partial \varphi}{\partial n^*}$:

$$\varphi_{n^*}(s,t) := n_2(s,t)\frac{\partial \varphi(s,t)}{\partial s} + n_1(s,t)\frac{\partial \varphi(s,t)}{\partial t}$$

--- (2.4)

The functional that is to be minimized is

$$J := \iint_G \Phi(s,t,x(s,t),x_s(s,t),x_t(s,t),u(s,t))\,dA(s,t) +$$

$$+ \int_{(\gamma)} \Phi_1(s,t,x(s,t),x_\mu(s,t))\,d\mu(s,t) + \sum_{P\in\mathbf{P}} \Phi_0(P,x(P))$$

--- (2.5)

where $d\mu(s,t) \equiv \sqrt{ds^2 + dt^2}$ is the differential of arc length on (γ). It is clear that the cost functional (2.4) includes (1.2) as a particular case.

The control function u takes values in $IR^m$, is assumed to be piecewise continuous in G (i.e. there is a finite collection of open sets, such that the union of their closures equals G, and u is continuous on each of those open sets), and for simplicity we assume that the components of u satisfy $\underline{u}^j \le u^j \le \overline{u}^j$, $1 \le j \le m$. The set of admissible values of the control function will be denoted by **U**. Naturally, more general controls can be considered, for example bounded measurable controls taking values in a compact subset of $IR^m$; however, such generality would not contribute to the issues of this paper, which concern the effects of the curvilinear boundary on the Hamiltonian equations.

3. Hamiltonian equations.

We derive the Hamiltonian equations arising out of the first variation of the functional (2.4) subject to the constraint of the state dynamics (2.1). By standard variational methods, the variation of the state, $\delta x(s,t)$, induced by an admissible variation $\delta u(s,t)$ of the control function, satisfies

$$\frac{\partial^2 \delta x^i(s,t)}{\partial s \partial t} = f^i_{x^j}(s,t,...)\delta x^j(s,t) + f^i_{p^j}(s,t,...)\delta x^j_s + f^i_{q^j}(s,t,...)\delta x^j_t(s,t) + $$
$$+ f^j_{u^k}(s,t,...)\delta u^k(s,t);\ \delta x^i(s,0) = \delta x^i(0,t) = 0$$

--- (3.1)

The variation of the cost functional J is

$$\delta J(u,\delta u) = \iint_G \{\Phi_{x^i}(s,t,...)\delta x^i(s,t) + \Phi_{p^i}(s,t,...)\delta x^i_s(s,t) + \Phi_{q^i}(s,t,...)\delta x^i_t(s,t) +$$
$$+ \Phi_{u^k}(s,t,...)\delta u^k(s,t)\} dA(s,t) +$$
$$+ \int_{(\gamma)} \{\Phi_{1,x^i}(s,t,...)\delta x^i(s,t) + \Phi_{1,\eta^i}(s,t,...)\delta x^i_\mu(s,t)\} d\mu(s,t) +$$
$$+ \sum_{r=1}^{N} \Phi_{0,x^i}(P_r, x(P_r))\delta x^i(P_r)$$

--- (3.2)

The variation $\delta x$ can be expressed in terms of the matrix-valued Riemann function $R^i_j(s,t,\sigma,\tau)$ associated with the linear Goursat-Darboux problem (3.1). The theory of the scalar Riemann function is explained in detail in [CC], and also it can be found in many standard references on partial differential equations. The proofs of the properties of the matrix-valued Riemann function are analogous to the proofs of [CC], and for this reason we shall only mention, without proofs, those of the properties of $[R^i_j]$ that are needed in this paper.
The matrix-valued Riemann function satisfies

$$\frac{\partial^2}{\partial s \partial t} R^i_j(s,t,\sigma,\tau) = f^i_{x^k}(s,t,...) R^k_j(s,t,\sigma,\tau) +$$

$$+ f^i_{p^k}(s,t,...) \frac{\partial}{\partial s} R^k_j(s,t,\sigma,\tau) + f^i_{q^k}(s,t,...) \frac{\partial}{\partial t} R^k_j(s,t,\sigma,\tau), \text{ for } s > \sigma \text{ and } t > \tau;$$

$$\frac{\partial}{\partial s} R^i_j(s,\tau,\sigma,\tau) = f^i_{q^k}(s,\tau,...) R^k_j(s,\tau,\sigma,\tau);$$

$$\frac{\partial}{\partial t} R^k_j(\sigma,t,\sigma,\tau) = f^i_{p^k}(\sigma,t,...) R^k_j(\sigma,t,\sigma,\tau);$$

$$R^i_j(\sigma,\tau,\sigma,\tau) = \delta^i_j \text{ (= Kronecker's delta)}$$

--- (3.3)

For later use, we also write down the adjoint problem to (3.3) that the matrix-valued Riemann function also solves:

$$\frac{\partial^2}{\partial s \partial t} R^i_j(\sigma,\tau,s,t) = \frac{\partial f^k(s,t,...)}{\partial x^j} R^i_k(\sigma,\tau,s,t) -$$

$$- \frac{D}{Ds}\left[\frac{\partial f^k(s,t,...)}{\partial p^j} R^i_k(\sigma,\tau,s,t)\right] - \frac{D}{Dt}\left[\frac{\partial f^k(s,t,...)}{\partial q^j} R^i_k(\sigma,\tau,s,t)\right],$$

for $s < \sigma$ and $t < \tau$;

$$\frac{\partial}{\partial t} R^i_j(\sigma,\tau,\sigma,t) = -\left[\frac{\partial f^k(\sigma,t,...)}{\partial p^j} R^i_k(\sigma,\tau,\sigma,t)\right];$$

$$\frac{\partial}{\partial s} R^i_j(\sigma,\tau,s,\tau) = -\left[\frac{\partial f^k(s,\tau,...)}{\partial q^j} R^i_k(\sigma,\tau,s,\tau)\right];$$

$$R^i_j(\sigma,\tau,\sigma,\tau) = \delta^i_j$$

--- (3.4)

The solution $\delta x(s,t)$ of (3.1) can be expressed as

$$\delta x^i(s,t) = \iint_{W(s,t)} R^i_j(s,t,\sigma,\tau) f^j_{u^k}(\sigma,\tau,...) \delta u^k(\sigma,\tau)\, dA(\sigma,\tau)$$

--- (3.5)

Thus, the variation of the cost functional J is

$$\delta J = \iint_G \iint_{W(s,t)} \{\Phi_{x^i}(s,t,...) - \frac{D}{Ds}\Phi_{p^i}(s,t,...) - \frac{D}{Dt}\Phi_{q^i}(s,t,...)\} \cdot$$
$$\cdot R^i_j(s,t,\sigma,\tau) f^j_{u^k}(\sigma,\tau,...)\delta u^k(\sigma,\tau)\, dA(\sigma,\tau)\, dA(s,t) +$$
$$+ \iint_G \Phi_{u^k}(s,t,...)\delta u^k(s,t)\, dA(s,t) +$$
$$+ \int_{(\gamma)} \iint_{W(s,t)} \{n_1(s,t)\Phi_{p^i}(s,t,...) + n_2(s,t)\Phi_{q^i}(s,t,...) + \Phi_{1,x^i}(s,t,...) -$$
$$- \frac{D}{D\mu}\Phi_{1,\eta^i}(s,t,...)\} R^i_j(s,t,\sigma,\tau) f^j_{u^k}(\sigma,\tau,...)\delta u^k(\sigma,\tau)\, dA(\sigma,\tau)\, d\mu(s,t) +$$
$$+ \sum_{r=1}^{N} \{\Phi_{0,x^i}(P_r, x(P_r)) - V_r(\Phi_{1,\eta^i})\} \iint_{W(P_r)} R^i_j(P_r,\sigma,\tau) f^j_{u^k}(\sigma,\tau,...)\delta u^k(\sigma,\tau)\, dA(\sigma,\tau)$$

--- (3.6)

By using the integration formulae of appendix A, we find

$$\delta J = \iint_G \iint_{E(s,t)} \{\Phi_{x^i}(\sigma,\tau,...) - \frac{D}{Ds}\Phi_{p^i}(\sigma,\tau,...) - \frac{D}{Dt}\Phi_{q^i}(\sigma,\tau,...)\} \cdot$$
$$\cdot R^i_j(\sigma,\tau,s,t) f^j_{u^k}(s,t,...)\delta u^k(s,t)\, dA(\sigma,\tau)\, dA(s,t) +$$
$$+ \iint_G \Phi_{u^k}(s,t,...)\delta u^k(s,t)\, dA(s,t) +$$
$$+ \iint_G \int_{\gamma(s,t)} \{n_1(\sigma,\tau)\Phi_{p^i}(\sigma,\tau,...) + n_2(\sigma,\tau)\Phi_{q^i}(\sigma,\tau,...) + \Phi_{1,x^i}(\sigma,\tau,...) -$$
$$- \frac{D}{D\mu}\Phi_{1,\eta^i}(\sigma,\tau,...)\} R^i_j(\sigma,\tau,s,t) f^j_{u^k}(s,t,...)\delta u^k(s,t)\, d\mu(\sigma,\tau)\, dA(s,t) +$$
$$+ \sum_{r=1}^{N} \{\Phi_{0,x^i}(P_r, x(P_r)) - V_r(\Phi_{1,\eta^i})\} \iint_G R^i_j(P_r,s,t) f^j_{u^k}(s,t,...)\delta u^k(s,t)\, dA(s,t)$$

--- (3.7)

The replacement of $\iint\limits_{W(P_r)} R^i_j(P_r,\sigma,\tau) f^j_{u^k}(\sigma,\tau,...)\delta u^k(\sigma,\tau)\, dA(\sigma,\tau)$ by $\iint\limits_G R^i_j(P_r,s,t) f^j_{u^k}(s,t,...)\delta u^k(s,t)\, dA(s,t)$, in the last term in (3.7) above, is justified by the fact that $R^i_j(P_r,s,t)=0$ whenever $(s,t)\notin W(P_r)$. For additional emphasis, we may write this term as $\iint\limits_G R^i_j(P_r,s,t) f^j_{u^k}(s,t,...)\chi((s,t)\in W(P_r))\delta u^k(s,t)\, dA(s,t)$ where $\chi$ stands for the truth function of a logical statement, $\chi(S) := \begin{cases} 1, & \text{if } S \text{ is true} \\ 0, & \text{if } S \text{ is false} \end{cases}$.

We define the co-state $\psi(s, t)$ by

$$
\begin{aligned}
\psi_j(s,t) :=& \iint\limits_{E(s,t)} \{\Phi_{x^i}(\sigma,\tau,...) - \frac{D}{Ds}\Phi_{p^i}(\sigma,\tau,...) - \frac{D}{Dt}\Phi_{q^i}(\sigma,\tau,...)\}\cdot \\
&\cdot R^i_j(\sigma,\tau,s,t)\, dA(\sigma,\tau) + \\
&+ \int\limits_{\gamma(s,t)} \{n_1(\sigma,\tau)\Phi_{p^i}(\sigma,\tau,...) + n_2(\sigma,\tau)\Phi_{q^i}(\sigma,\tau,...) + \Phi_{1,x^i}(\sigma,\tau,...) - \\
&- \frac{D}{D\mu}\Phi_{1,\eta^i}(\sigma,\tau,...)\} R^i_j(\sigma,\tau,s,t)\, d\mu(\sigma,\tau) + \\
&+ \sum_{r=1}^{N} \{\Phi_{0,x^i}(P_r, x(P_r)) - V_r(\Phi_{1,\eta^i})\} R^i_j(P_r,s,t)\chi((s,t)\in W(P_r))
\end{aligned}
$$

--- (3.8)

Also, we define the Hamiltonian H by

$$
H(s,t,x,p,q,\psi,u) := \Phi(s,t,x,p,q,u) + \psi_j f^j(s,t,x,p,q,u)
$$

--- (3.9)

Then the variation of J takes the form

$$
\delta J = \iint\limits_G \frac{\partial H(s,t,x(s,t),x_s(s,t),x_t(s,t),\psi(s,t),u(s,t))}{\partial u^k} \delta u^k(s,t)\, dA(s,t)
$$

--- (3.10)

We shall prove

<u>Theorem 3.1.</u> Assuming adequate continuous differentiability of all functions involved, the co-state $\psi$ defined in (3.8) satisfies the following hyperbolic equation, in the regular part of the domain G:

$$\frac{\partial^2 \psi_i(s,t)}{\partial s \partial t} = \frac{\partial H(s,t,x(s,t),x_s(s,t),x_t(s,t),\psi(s,t),u(s,t))}{\partial x^i} -$$
$$- \frac{D}{Ds}\frac{\partial H(s,t,x(s,t),x_s(s,t),x_t(s,t),\psi(s,t),u(s,t))}{\partial p^i} -$$
$$- \frac{D}{Dt}\frac{\partial H(s,t,x(s,t),x_s(s,t),x_t(s,t),\psi(s,t),u(s,t))}{\partial q^i}$$

--- (3.11)

<u>Proof:</u> We shall use the adjoint equations (3.4) for the Riemann function and the differentiation formulae of appendix A. In order to simplify the notation, we set

$$F_i(\sigma,\tau,...) := \Phi_{x^i}(\sigma,\tau,...) - \frac{D}{D\sigma}\Phi_{p^i}(\sigma,\tau,...) - \frac{D}{D\tau}\Phi_{q^i}(\sigma,\tau,...);$$
$$F_{1,i}(\sigma,\tau,...) := n_1(\sigma,\tau)\Phi_{p^i}(\sigma,\tau,...) + n_2(\sigma,\tau)\Phi_{q^i}(\sigma,\tau,...) + \Phi_{1,x^i}(\sigma,\tau,...) -$$
$$- \frac{D}{D\mu}\Phi_{1,\eta^i}(\sigma,\tau,...);$$
$$F_{0,i}(P_r, x(P_r)) := \Phi_{0,x^i}(P_r, x(P_r)) - V_r(\Phi_{1,\eta^i})$$

--- (3.12)

Then (3.8) can be written as

$$\psi_j(s,t) = \iint_{E(s,t)} F_i(\sigma,\tau,...) R^i_j(\sigma,\tau,s,t)\, dA(\sigma,\tau) +$$
$$+ \int_{\gamma(s,t)} F_{1,i}(\sigma,\tau,...) R^i_j(\sigma,\tau,s,t)\, d\mu(\sigma,\tau) + \sum_{r=1}^{N} F_{0,i}(P_r, x(P_r)) R^i_j(P_r,s,t)$$

--- (3.13)

According to the differentiation formulae that are proved in appendix A, we have

$$\frac{\partial \psi_j(s,t)}{\partial s} = -\int_{P_{st}}^{B_s} F_i(s,\tau,...)R_j^i(s,\tau,s,t)\,d\tau +$$

$$+ \int_{E(s,t)} F_1(\sigma,\tau,...)\frac{\partial}{\partial s}R_j^i(\sigma,\tau,s,t)\,dA(\sigma,\tau) - F_{1,i}(B_s,...)R_j^i(B_s,s,t)\frac{1}{n_2(B_s)} +$$

$$+ \int_{(\gamma)} F_{1,i}(\sigma,\tau,...)\frac{\partial}{\partial s}R_j^i(\sigma,\tau,s,t)\,d\mu(\sigma,\tau) +$$

$$+ \sum_{r=1}^{N} F_{0,i}(P_r,x(P_r))\chi((s,t)\in W(P_r))\frac{\partial}{\partial s}R_j^i(P_r,s,t)$$

--- (3.14)

$$\frac{\partial \psi_j(s,t)}{\partial t} = -\int_{P_{st}}^{A_t} F_i(s,\tau,...)R_j^i(s,\tau,s,t)\,d\sigma +$$

$$+ \int_{E(s,t)} F_1(\sigma,\tau,...)\frac{\partial}{\partial t}R_j^i(\sigma,\tau,s,t)\,dA(\sigma,\tau) - F_{1,i}(A_t,...)R_j^i(A_t,s,t)\frac{1}{n_1(A_t)} +$$

$$+ \int_{(\gamma)} F_{1,i}(\sigma,\tau,...)\frac{\partial}{\partial t}R_j^i(\sigma,\tau,s,t)\,d\mu(\sigma,\tau) +$$

$$+ \sum_{r=1}^{N} F_{0,i}(P_r,x(P_r))\chi((s,t)\in W(P_r))\frac{\partial}{\partial t}R_j^i(P_r,s,t)$$

--- (3.15)

$$\frac{\partial^2\psi_j(s,t)}{\partial s\partial t} = -\int_{P_{st}}^{A_t} F_i(s,\tau,...)R^i_j(s,\tau,s,t)\,d\sigma - \int_{P_{st}}^{B_s} F_i(s,\tau,...)R^i_j(s,\tau,s,t)\,d\tau +$$

$$+ \int_{E(s,t)} F_1(\sigma,\tau,...)\frac{\partial^2}{\partial s\partial t}R^i_j(\sigma,\tau,s,t)\,dA(\sigma,\tau) -$$

$$- F_{1,i}(A_t,...)R^i_j(A_t,s,t)\frac{1}{n_1(A_t)} - F_{1,i}(B_s,...)R^i_j(B_s,s,t)\frac{1}{n_2(B_s)} +$$

$$+ \int_{(\gamma)} F_{1,i}(\sigma,\tau,...)\frac{\partial^2}{\partial s\partial t}R^i_j(\sigma,\tau,s,t)\,d\mu(\sigma,\tau) +$$

$$+ \sum_{r=1}^{N} F_{0,i}(P_r,x(P_r))\chi((s,t)\in W(P_r))\frac{\partial^2}{\partial s\partial t}R^i_j(P_r,s,t)$$

--- (3.16)

We note that, in the interior of $G_R$, the function $\chi((s,t)\in W(P_r))$ remains constant. For those terms of (3.16) that contain $\frac{\partial^2}{\partial s\partial t}R^i_j(\sigma,\tau,s,t)$, we use the expanded form of the first equation in (3.4):

$$\frac{\partial^2}{\partial s\partial t}R^i_j(\sigma,\tau,s,t) = \frac{\partial f^k(s,t,...)}{\partial x^j}R^i_k(\sigma,\tau,s,t) - \frac{\partial f^k(s,t,...)}{\partial p^i}\frac{\partial}{\partial s}R^i_j(\sigma,\tau,s,t) -$$

$$- \frac{\partial f^k(s,t,...)}{\partial q^i}\frac{\partial}{\partial t}R^i_j(\sigma,\tau,s,t) - R^i_k(\sigma,\tau,s,t)\frac{D}{Ds}\frac{\partial f^k(s,t,...)}{\partial p^j} - R^i_k(\sigma,\tau,s,t)\frac{D}{Dt}\frac{\partial f^k(s,t,...)}{\partial q^j}$$

--- (3.17)

The expanded form of the wanted equation (3.11) is

$$\frac{\partial^2\psi_j(s,t)}{\partial s\partial t} = \frac{\partial\Phi(s,t,...)}{\partial x^j} - \frac{D}{Ds}\frac{\partial\Phi(s,t,...)}{\partial p^j} - \frac{D}{Dt}\frac{\partial\Phi(s,t,...)}{\partial q^j} +$$

$$+ \psi_i(s,t)\left[\frac{\partial f^i(s,t,...)}{\partial x^j} - \frac{D}{Ds}\frac{\partial f^i(s,t,...)}{\partial p^j} - \frac{D}{Dt}\frac{\partial f^i(s,t,...)}{\partial q^j}\right] -$$

$$- \frac{\partial\psi_i(s,t)}{\partial s}\frac{\partial f^i(s,t,...)}{\partial p^j} - \frac{\partial\psi_i(s,t)}{\partial t}\frac{\partial f^i(s,t,...)}{\partial q^j}$$

--- (3.18)

We examine the terms that arise out of the right-hand side of (3.16), and, using (3.17) and the remaining equations out of (3.4), we show that these terms match the terms on the right-hand side of the wanted equation (3.18).

Each term containing $\frac{\partial^2 R^i_j}{\partial s \partial t}$ on the right-hand side of (3.16), by invoking (3.17), is matched with the corresponding terms containing $\frac{\partial R^i_k}{\partial s}, \frac{\partial R^i_k}{\partial t}$ multiplied by $\frac{\partial f^k}{\partial p^j}, \frac{\partial f^k}{\partial q^j}$, arising from the terms on the right-hand side of (3.18) that contain $\frac{\partial \psi_i}{\partial s}, \frac{\partial \psi_i}{\partial t}$, after those terms have been expressed by utilizing the second and third equations out of the set (3.4), and with similar terms containing $R^i_k$ multiplied by $\frac{D}{Ds}\frac{\partial f^k}{\partial p^j}, \frac{D}{Dt}\frac{\partial f^k}{\partial q^j}$ arising from the terms containing $\psi_i$ on the right-hand side of (3.16). This matching can be easily seen, for example, for the terms consisting of double integrals over E(s, t) ; the proof is similar for the other terms that contain $\frac{\partial^2 R^i_j}{\partial s \partial t}$ on the right-hand side of (3.16). Now, by the third equation in (3.4), the term $\frac{\partial}{\partial s} R^i_j(\sigma, \tau, s, t)$ equals $-\frac{\partial f^k(s, t, ...)}{\partial q^j} R^i_k(\sigma, \tau, s, t)$, and this is matched with the corresponding term in the expansion of $-\frac{\partial \psi_i}{\partial t}\frac{\partial f^i}{\partial q^j}$. This comparison is typified by the terms containing the integrals $\int_{P_{st}}^{A_t}$ in (3.16) after $\psi_i$ and its first derivatives (from (3.14) and (3.15)) have been substituted into (3.16); the comparison for the remaining terms is similar. Of course, the comparison for the terms that contain $\frac{\partial R^i_j}{\partial t}$ is parallel to the comparison, explained above, for the terms that contain $\frac{\partial R^i_j}{\partial s}$. ///

4. Side conditions for the Hamiltonian equations.

The solution of the Hamiltonian equations, which were derived in the previous section, requires side conditions that make the corresponding Goursat problems (Hamiltonian equations plus side conditions) well-posed problems. These side conditions are not merely boundary conditions but rather they are more complicated conditions, and they constitute the main difference between the case of rectangular domains and the case of general domains with curvilinear boundary of the type described in this paper.

We shall need some definitions and notation pertaining to the geometry of $(\gamma)$.

Definition 4.1. A connected part of $(\gamma)$ will be called flat if it is a straight line segment parallel to either the s-axis or the t-axis; a connected part of $(\gamma)$ will be called oblique if it does not contain any straight line segment that is parallel to either the s-axis or the t-axis. ///

Theorem 4.1. On the relative interior (in the relative topology of $(\gamma)$ ) of a flat part of $(\gamma)$, at points (s, t) that do not belong to **P**, the co-state $\psi$ satisfies

$$\left[(n_1(s,t))^2 - (n_2(s,t))^2\right]\frac{\partial\psi_j(s,t)}{\partial\mu} + n_1(s,t)\frac{\partial H(s,t,...)}{\partial p^j} + n_2(s,t)\frac{\partial H(s,t,...)}{\partial q^j} +$$

$$+\frac{\partial\Phi_1(s,t,...)}{\partial x^j} - \frac{D}{D\mu}\frac{\partial\Phi_1(s,t,...)}{\partial\eta^j} = 0$$

--- (4.1)

If P is a point in the relative interior of the oblique part of $(\gamma)$ and $P \notin \mathbf{P}$, then

$$\lim_{(s,t)\to P}\psi_j(s,t) = 0;$$

$$\lim_{(s,t)\to P}\left[n_2(P)\frac{\partial\psi_j(s,t)}{\partial s}\right] = \lim_{(s,t)\to P}\left[n_1(P)\frac{\partial\psi_j(s,t)}{\partial t}\right] = -F_{1,j}(P,...)$$

--- (4.2)

Proof: For a point (s, t) in the relative interior of a flat part of $(\gamma)$ with $(s,t) \notin \mathbf{P}$, we have, in some relative neighbourhood of (s, t), either s=const. or t=const., and each of the functions $\chi((\sigma,\tau) \in W(P_r))$ remains constant for $(\sigma, \tau)$ in a sufficiently small relative neighbourhood of (s, t). Clearly, it suffices to prove (4.1) in one of the cases s=const. or t=const., for example in the case $s = s_1 = \text{const.}$ In that case, we have $B_s \equiv (s_1, t_1)$, for some value $t_1$ such that $\{(s_1,t): t_0 \le t \le t_1\}$ is a maximal straight line segment containing the point (s, t), and $B_s$ is a member of **P**, according to our definitions and assumptions in section 2.

There can be other elements, say $P_{r'}$, of **P** in the segment of (γ) that connects $(s_1,t)$ and $B_s$. On that segment, we have $d\mu \equiv dt$, $n_1(s,t) \equiv 1$, $n_2(s,t) \equiv 0$. As $(s,t) \to (s_1,t)$, with (s, t) in the interior of the regular part of G, the double integral over E(s, t) goes to 0, and we have, according to the definition of ψ,

$$\psi_j(s_1,t) = \int_t^{t_1} F_i(s_1,\tau,...)R^i_j(s_1,\tau,s_1,t)\,d\tau + F_{0,i}((s_1,t_1),x(s_1,t_1))R^i_j(s_1,t_1,s_1,t) +$$

$$+\sum_{r'} F_{0,i}((s_1,t(P_{r'})),x(s_1,t(P_{r'})))R^i_j(s_1,t(P_{r'}),s_1,t)$$

--- (4.3)

form which

$$\frac{\partial\psi_j(s_1,t)}{\partial t} = -F_i(s_1,t,...)R^i_j(s_1,t,s_1,t) + \int_t^{t_1} F_i(s_1,t,...)\frac{\partial}{\partial t}R^i_j(s_1,\tau,s_1,t)\,d\tau +$$

$$+F_{0,i}((s_1,t_1),x(s_1,t_1))\frac{\partial}{\partial t}R^i_j(s_1,t_1,s_1,t) +$$

$$+\sum_{r'} F_{0,i}((s_1,t(P_{r'})),x(s_1,t(P_{r'})))\frac{\partial}{\partial t}R^i_j(s_1,t(P_{r'}),s_1,t)$$

--- (4.4)

By using (3.4), we can write (4.4) as

$$\frac{\partial\psi_j(s_1,t)}{\partial t} = -F_j(s_1,t,...) - \int_t^{t_1} F_i(s_1,t,...)\frac{\partial f^k(s_1,t,...)}{\partial p^j}R^i_k(s_1,\tau,s_1,t)\,d\tau +$$

$$+F_{0,i}((s_1,t_1),x(s_1,t_1))\frac{\partial f^k(s_1,t,...)}{\partial p^j}R^i_k(s_1,t_1,s_1,t) +$$

$$+\sum_{r'} F_{0,i}((s_1,t(P_{r'})),x(s_1,t(P_{r'})))\frac{\partial f^k(s_1,t,...)}{\partial p^j}R^i_k(s_1,t(P_{r'}),s_1,t) =$$

$$= -F_j(s_1,t,...) - \frac{\partial}{\partial p^j}\left[\psi_k(s,t)f^k(s_1,t,...)\right]$$

--- (4.5)

Now, (4.5) is tantamount to (4.1); this last assertion is a direct consequence of the definitions of H and $F_{1,i}$ in (3.9) and (3.12).

If $(s,t)$ is a point on the oblique part of ($\gamma$) and not on **P**, then, as $(s',t')\to(s,t)$ with $(s',t')\in \mathrm{int}(G_R)$, both the double integral over E(s', t') and the line integral over $\gamma$(s', t'), in the definition of $\psi$(s', t'), will go to 0; the discrete terms containing $F_{0,i}$ will also be zero if (s', t') is sufficiently close to (s, t) so that there is no point P in **P** such that W(P) contains (s', t'). This proves the first equation in (4.2).

To prove the second equation in (4.2), we use (3.14) and (3.15); as $(s,t)\to P\equiv(s(P),t(P))$, the only terms, out of the right-hand sides of (3.14) and (3.15), that have (possibly) nonzero limits are the terms

$-F_{1,i}(B_s,...)R^i_j(B_s,s,t)\dfrac{1}{n_2(B_s)}$, $-F_{1,i}(A_t,...)R^i_j(A_t,s,t)\dfrac{1}{n_1(A_t)}$ ; thus

$$\lim_{(s,t)\to P}\left[n_2(P)\frac{\partial\psi_j(s,t)}{\partial s}\right]=\lim_{(s,t)\to P}\left[-F_{1,i}(B_s,...)R^i_j(B_s,s,t)\frac{n_2(P)}{n_2(B_s)}\right]=-F_{1,j}(P,...) \qquad \text{--- (4.6)}$$

and analogously

$$\lim_{(s,t)\to P}\left[n_1(P)\frac{\partial\psi_j(s,t)}{\partial t}\right]=-F_{1,j}(P,...) \qquad \text{--- (4.7)}$$

This concludes the proof of (4.2). ///

The examination of the limiting values of $\psi$ at a point of **P** requires the following

Definition 4.2. For every point P of G (including points on $(\gamma_R)$ and points on **P**), we define the sets

$W^{III}(P):=\{(s,t)\in \mathrm{int}(G): s<s(P),\ t<t(P)\}$;

$W^{II}(P):=\{(s,t)\in \mathrm{int}(G): s<s(P),\ t>t(P)\}$;

$W^{IV}(P):=\{(s,t)\in \mathrm{int}(G): s>s(P),\ t<t(P)\}$

--- (4.8)

(The superscripts in (4.8) refer to the quadrants of the st- plane around the point P.)

For every point P on ($\gamma$), we define the following sets:

T(P) is the collection of all vertices, excluding P itself in case P is a vertex, subsequent to P relative to the counterclockwise orientation of $\partial G$, that lie on a straight line segment parallel to the t-axis, that is part of ($\gamma$) and passes through P;
S(P) is the collection of all vertices, excluding P itself in case P is a vertex, preceding P relative to the counterclockwise orientation of $\partial G$, that lie on a straight line segment parallel to the s-axis, that is part of ($\gamma$) and passes through P;
$V(P) := T(P) \cup S(P) \cup \{P\}$;
$L_t(P)$ is the maximal straight line segment contained in $(\gamma)$, excluding P itself, consisting of points subsequent to P relative to the counterclockwise orientation of $\partial G$, that is parallel to the t-axis and passes through P;
$L_s(P)$ is the maximal straight line segment contained in $(\gamma)$, excluding P itself, consisting of points preceding P relative to the counterclockwise orientation of $\partial G$, that is parallel to the s-axis and passes through P;
$L(P) := L_s(P) \cup L_t(P) \cup \{P\}$.

Clearly, in some cases, some of the sets defined here may be empty. ///

We have:

<u>Theorem 4.2.</u> The limiting values of $\psi$ as (s, t) approaches a vertex $P_k$ are determined by

$$\lim_{\substack{(s,t)\to P_k \\ (s,t)\in W^{III}(P_k)}} \psi_j(s,t) = \int_{L(P_k)} F_{1,i}(\sigma,\tau,...)R_j^i(\sigma,\tau,P_k)\,d\mu(\sigma,\tau) + $$

$$+ \sum_{P_{k'}\in V(P_k)} F_{0,i}(P_{k'},...)R_j^i(P_{k'},P_k);$$

$$\lim_{\substack{(s,t)\to P_k \\ (s,t)\in W^{II}(P_k)}} \psi_j(s,t) = \int_{L_t(P_k)} F_{1,i}(\sigma,\tau,...)R_j^i(\sigma,\tau,P_k)\,d\mu(\sigma,\tau) + $$

$$+ \sum_{P_{k'}\in T(P_k)} F_{0,i}(P_{k'},...)R_j^i(P_{k'},P_k);$$

$$\lim_{\substack{(s,t)\to P_k \\ (s,t)\in W^{IV}(P_k)}} \psi_j(s,t) = \int_{L_s(P_k)} F_{1,i}(\sigma,\tau,...)R_j^i(\sigma,\tau,P_k)\,d\mu(\sigma,\tau) + $$

$$+ \sum_{P_{k'}\in S(P_k)} F_{0,i}(P_{k'},...)R_j^i(P_{k'},P_k)$$

--- (4.9)

<u>Proof:</u> For (s, t) in $W^{III}(P_k)$ and sufficiently close to $P_k$, we have $E(s,t)\cap \mathbf{P} = L(P_k)$. As $(s,t)\to P_k$ with $(s,t)\in W^{III}(P_k)$, the arc $A_tB_s$ converges to $L(P_k)$ in the sense that the two points $A_t$, $B_s$ converge to the two extremities of $L(P_k)$. To show this, we observe that, if $L(P_k) = P_{k'}P_{k''}$ (i.e. the set $L(P_k)$ is equal to the connected arc of ($\gamma$) that has extremities $P_{k'}$ and $P_{k''}$), with $P_{k'}$ preceding $P_{k''}$, then the segment $P_{k''}P_k$ is a maximal straight-line segment through $P_k$ that is part of ($\gamma$) and is parallel to the s-axis, and $P_kP_{k'}$ is a maximal straight-line segment through $P_k$ that is part of ($\gamma$) and is parallel to the t-axis; this situation includes the possibility that one or both of the straight-line segments $P_{k''}P_k$, $P_kP_{k'}$ might degenerate into the singleton $\{P_k\}$. By definition, for $(s,t)\in W^{III}(P_k)$, we have $s < s(P_k)$ and $t < t(P_k)$, thus also $s < s(P_{k'}) = s(P_k)$ and $t < t(P_{k''}) = t(P_k)$, and consequently $E(s,t) \supseteq L(P_k)$. By the maximality of the segments $P_{k''}P_k$, $P_kP_{k'}$ , the arc $A_tP_{k''}$ is not parallel to the s-axis, and the arc $P_{k'}B_s$ is not parallel to the t-axis. Therefore, the point $B_s$ has the representation $B_s = (s, \vartheta_2(s))$ and consequently, as $s \to s(P_k) = s(P_{k''})$, we have $B_s \to P_{k''}$. By the same token, as $t \to t(P_k) = t(P_{k'})$, we have $A_t = (\vartheta_1(t), t) \to P_{k'}$. By taking limits in (3.8), which is the the definition of $\psi$, using the information above, we obtain the first equation in (4.9). The proofs of the remaining parts of (4.9) are similar. ///

Definition 4.3. We denote by $(S_j)$ the straight line segment through $P_j$ and parallel to the s-axis, by $(T_k)$ the straight line segment through $P_k$ and parallel to the t-axis, and by $V_{j,k}$ the intersection of $(S_j)$ and $(T_k)$ for k>j. We denote by $D_j$, $0 \le j \le N$, the (generally) curvilinear-triangular domain bounded by the sub-arc $P_jP_{j+1} \subseteq (\gamma)$ and the straight line segments $V_{j,j+1}P_j$, $V_{j,j+1}P_{j+1}$. (Clearly, $D_j$ degenerates into the straight-line segment $P_jP_{j+1}$ when that segment is a flat part of (γ), and therefore, in that case, $\operatorname{int} D_j = \varnothing$ .) We denote by $Q_{j,k}$, $0 \le j < k \le N$ , the rectangular domain with vertices $V_{j,k}$, $V_{j+1,k}$, $V_{j,k+1}$, $V_{j+1,k+1}$. The collection of domains $\{D_j, Q_{j,k}\}$ is partitioned into zones as follows: zone $Z_0$ comprises the triangular domains $D_j$; for each integer p, $1 \le p \le N$, zone $Z_p$ comprises all rectangular domains $Q_{j,k}$ with $k - j = p$. ///

Definition 4.4. The jump $\Omega_{(T)}\varphi$ of a function φ across a straight line segment (T) that is parallel to the t-axis is defined as $\Omega_{(T)}\varphi(s,t) := \varphi^-_{(T)}(s,t) - \varphi^+_{(T)}(s,t)$ for $(s,t) \in (T)$ , where $\varphi^-_{(T)}(s,t) := \lim\limits_{\sigma \to s^-} \varphi(\sigma,t)$, $\varphi^+_{(T)}(s,t) := \lim\limits_{\sigma \to s^+} \varphi(\sigma,t)$. Similarly, the jump $\Omega_{(S)}\varphi$ of a function φ across a straight line segment (S) that is parallel to the s-axis is defined as $\Omega_{(S)}\varphi(s,t) := \varphi^-_{(S)}(s,t) - \varphi^+_{(S)}(s,t)$, where $\varphi^-_{(S)}(s,t) := \lim\limits_{\tau \to t^-} \varphi(s,\tau)$, $\varphi^+_{(S)}(s,t) := \lim\limits_{\tau \to t^+} \varphi(s,\tau)$, for $(s,t) \in (S)$. Of course, these definitions are conditional on the existence of the indicated limits. ///

Theorem 4.3. The jumps of the co-state ψ across the lines $(S_r)$, $(T_r)$ are given by

$$\Omega_{(S_r)}\psi_j(s,t) = \int\limits_{L_s(P_r)} F_{1,i}(\sigma,\tau,\ldots)R^i_j(\sigma,\tau,s,t)\,d\mu(\sigma,\tau) + \sum_{r':P_{r'}\in S(P_r)} F_{0,i}(P_{r'},\ldots)R^i_j(P_{r'},s,t)$$

for $(s,t) \in \operatorname{int} G$ with $t \ne t(P_k)$ $\forall k$;

$$\Omega_{(T_r)}\psi_j(s,t) = \int\limits_{L_t(P_r)} F_{1,i}(\sigma,\tau,\ldots)R^i_j(\sigma,\tau,s,t)\,d\mu(\sigma,\tau) + \sum_{r'':P_{r''}\in S(P_r)} F_{0,i}(P_{r''},\ldots)R^i_j(P_{r''},s,t)$$

for $(s,t) \in \operatorname{int} G$ with $s \ne s(P_k)$ $\forall k$

--- (4.10)

Proof: As in the proof of theorem 4.2, the contribution of the double integrals (integrals over E(s, t)) to the calculation of the jump of ψ is zero. We prove the first equation in (4.10), since the two equations in (4.10) are symmetric with respect to an interchange of the variables s and t. We consider two straight lines $S_{r,\varepsilon}$, $S_{r,-\varepsilon}$ consisting of points of the form $\{(s(P_r),t'): t'=t+\varepsilon,\ (s(P_r),t)\in S_r\}$, $\{(s(P_r),t''): t''=t-\varepsilon,\ (s(P_r),t)\in S_r\}$, respectively. For ε>0 and sufficiently small, the points of **P** that are included in $\gamma(s(P_r),t-\varepsilon)$ but not in $\gamma(s(P_r),t+\varepsilon)$ are precisely the points in $L_s(P_r)$. If the vertices $P_{r'}$, $P_{r''}$ are defined in the same way as the points $P_{k'}$, $P_{k''}$, respectively, in the proof of theorem 4.2, then the set-theoretic difference $\gamma(s(P_r),t-\varepsilon)\setminus\gamma(s(P_r),t+\varepsilon)$ converges, as $\varepsilon\to 0^+$, to the arc $P_{r'}P_r$ which coincides with the set $S(P_r)$. Therefore, we have

$$
\begin{aligned}
&\Omega_{(S_r)}\psi_j(s(P_r),t) = \\
&= \lim_{\varepsilon\to 0^+}\Bigg[\iint_{E(s(P_r),t-\varepsilon)} F_i(\sigma,\tau,\ldots)R^i_j(\sigma,\tau,s(P_r),t-\varepsilon)\,dA(\sigma,\tau) - \\
&- \iint_{E(s(P_r),t+\varepsilon)} F_i(\sigma,\tau,\ldots)R^i_j(\sigma,\tau,s(P_r),t+\varepsilon)\,dA(\sigma,\tau) + \\
&+ \int_{\gamma(s(P_r),t-\varepsilon)} F_{1,i}(\sigma,\tau,\ldots)R^i_j(\sigma,\tau,s(P_r),t-\varepsilon)\,d\mu(\sigma,\tau) - \\
&- \int_{\gamma(s(P_r),t+\varepsilon)} F_{1,i}(\sigma,\tau,\ldots)R^i_j(\sigma,\tau,s(P_r),t+\varepsilon)\,d\mu(\sigma,\tau) + \\
&+ \sum_{r':\,P_{r'}\in S(P_r)} F_{0,i}(P_{r'},\ldots)R^i_j(P_{r'},s(P_r),t-\varepsilon) + \\
&+ \sum_{k:\,P_k\in W(s(P_r),t+\varepsilon)} F_{0,i}(P_k,\ldots)\Big(R^i_j(P_k,s(P_r),t-\varepsilon) - R^i_j(P_k,s(P_r),t+\varepsilon)\Big)\Bigg] = \\
&= \int_{L_s(P_r)} F_{1,i}(\sigma,\tau,\ldots)R^i_j(\sigma,\tau,s(P_r),t)\,d\mu(\sigma,\tau) + \sum_{r':\,P_{r'}\in S(P_r)} F_{0,i}(P_{r'},\ldots)R^i_j(P_{r'},s(P_r),t)
\end{aligned}
$$

--- (4.11)

Theorem 4.3 is proved. ///

The side conditions of theorems 4.2 and 4.3 can be written in Hamiltonian form. The side conditions depend on the Riemann function through terms of the form $R^i_j(s,\tau,s,t)$, $R^i_j(\sigma,t,s,t)$. We set

$$\rho^i_{1,j}(s,t,\tau) := R^i_j(s,\tau,s,t),\ \rho^i_{2,j}(s,t,\sigma) := R^i_j(\sigma,t,s,t)$$

--- (4.12)

The quantities $\rho^i_{\alpha,j}$, $\alpha = 1,2$, will play the role of additional new co-states. According to (3.2), we have

$$\frac{\partial}{\partial t}\rho^i_{1,j}(s,t,\tau) = -\frac{\partial f^k(s,t,...)}{\partial p^j}\rho^i_{1,k}(s,t,\tau);\ \rho^i_{1,j}(s,\tau,\tau) = \delta^i_j;$$
$$\frac{\partial}{\partial s}\rho^i_{2,j}(s,t,\sigma) = -\frac{\partial f^k(s,t,...)}{\partial q^j}\rho^i_{2,k}(s,t,\sigma);\ \rho^i_{2,j}(\sigma,t,\sigma) = \delta^i_j$$

--- (4.13)

We define two new families of Hamiltonians $h^i_1$, $h^i_2$, $1 \le i \le n$ as follows:

$$h^i_1(s,t,\tau,x,p,q,\rho_1,u) := f^k(s,t,x,p,q,u)\rho^i_{1,k};$$
$$h^i_2(s,t,\tau,x,p,q,\rho_2,u) := f^k(s,t,x,p,q,u)\rho^i_{2,k};$$

--- (4.14)

Then (4.13) can be written in Hamiltonian form as

$$\frac{\partial}{\partial t}\rho^i_{1,j}(s,t,\tau) = -\frac{\partial h^i_1(s,t,\tau,x(s,t),x_s(s,t),x_t(s,t),\rho_1(s,t,\tau),u(s,t))}{\partial p^j};$$
$$\frac{\partial}{\partial s}\rho^i_{2,j}(s,t,\sigma) = -\frac{\partial h^i_2(s,t,\tau,x(s,t),x_s(s,t),x_t(s,t),\rho_2(s,t,\tau),u(s,t))}{\partial p^j};$$
$$\rho^i_{1,j}(s,\tau,\tau) = \rho^i_{2,j}(\sigma,t,\sigma) = \delta^i_j$$

--- (4.15)

The conditions of theorems 4.2 and 4.3 become

$$\lim_{\substack{(s,t)\to P_k \\ (s,t)\in W^{III}(P_k)}} \psi_j(s,t) = \int_{s(P_k)}^{s(P_{k'})} F_{1,i}(\sigma, t(P_k),\ldots)\rho^i_{2,j}(s(P_k), t(P_k), \sigma)\, d\sigma +$$

$$+ \int_{t(P_k)}^{t(P_{k''})} F_{1,i}(s(P_k), \tau,\ldots)\rho^i_{1,j}(s(P_k), t(P_k), \tau)\, d\tau +$$

$$+ \sum_{k' \le r < k} F_{0,i}(s(P_r), t(P_k),\ldots)\rho^i_{2,j}(s(P_k), t(P_k), s(P_r)) +$$

$$+ \sum_{k < r \le k''} F_{0,i}(s(P_k), t(P_r),\ldots)\rho^i_{1,j}(s(P_k), t(P_k), t(P_r)) + F_{0,j}(s(P_k), t(P_k),\ldots);$$

$$\lim_{\substack{(s,t)\to P_k \\ (s,t)\in W^{II}(P_k)}} \psi_j(s,t) = \int_{s(P_k)}^{s(P_{k'})} F_{1,i}(\sigma, t(P_k),\ldots)\rho^i_{2,j}(s(P_k), t(P_k), \sigma)\, d\sigma +$$

$$+ \sum_{k' \le r < k} F_{0,i}(s(P_r), t(P_k),\ldots)\rho^i_{2,j}(s(P_k), t(P_k), s(P_r));$$

$$\lim_{\substack{(s,t)\to P_k \\ (s,t)\in W^{IV}(P_k)}} \psi_j(s,t) = \int_{t(P_k)}^{t(P_{k''})} F_{1,i}(s(P_k), \tau,\ldots)\rho^i_{1,j}(s(P_k), t(P_k), \tau)\, d\tau +$$

$$+ \sum_{k < r \le k''} F_{0,i}(s(P_k), t(P_r),\ldots)\rho^i_{1,j}(s(P_k), t(P_k), t(P_r))$$

--- (4.16)

$$\Omega_{(S_k)}\psi_j(s,t) = \int_{s(P_{k'})}^{s(P_k)} F_{1,i}(\sigma, t(P_k),\ldots)\rho^i_{2,j}(s(P_k), t(P_k), \sigma)\, d\sigma +$$

$$+ \sum_{k' \le r < k} F_{0,i}(s(P_r), t(P_k),\ldots)\rho^i_{2,j}(s(P_k), t(P_k), s(P_r));$$

$$\Omega_{(T_k)}\psi_j(s,t) = \int_{t(P_{k'})}^{t(P_k)} F_{1,i}(s(P_k), \tau,\ldots)\rho^i_{1,j}(s(P_k), t(P_k), \tau)\, d\tau +$$

$$+ \sum_{k < r \le k''} F_{0,i}(s(P_k), t(P_r),\ldots)\rho^i_{1,j}(s(P_k), t(P_k), t(P_r))$$

--- (4.17)

The significance of the side conditions in theorems 4.1, 4.2, and 4.3 is shown in the following:

Theorem 4.4. We assume that the right-hand side of the Hamiltonian equations (3.12) is Lipschitz continuous with respect to $\left(\psi, \frac{\partial\psi}{\partial s}, \frac{\partial\psi}{\partial t}\right)$, uniformly in all the other variables. Then the side conditions for $\psi$, established in theorems 4.1, 4.2, and 4.3, are sufficient for the unique solvability of the Hamiltonian equations (3.12).

Proof: For each non-degenerate triangular domain $D_j$, the Hamiltonian equations (3.12) are written, in integral form (after repeated application of the Gauss-Green theorem), as

$$\psi_i(s,t) = \frac{1}{2}\left(\psi_i(B_s) + \psi_i(A_t)\right) + \iint\limits_{E(s,t)} g_i(\sigma,\tau,\psi(\sigma,\tau), \frac{\partial\psi(\sigma,\tau)}{\partial\tau}, \frac{\partial\psi(\sigma,\tau)}{\partial\sigma}) dA(\sigma,\tau) -$$
$$-\frac{1}{2}\int_{A_t}^{B_s} \psi_i^{n^*}(\sigma,\tau)\, d\mu(\sigma,\tau)$$

--- (4.18)

where $g_i$ denotes the right-hand side of (3.12), i.e.

$$g_i(s,t,\psi,\psi_s,\psi_t) := \frac{\partial H(s,t,x(s,t),x_s(s,t),x_t(s,t),\psi(s,t),u(s,t))}{\partial x^i} -$$
$$-\frac{D}{Ds}\frac{\partial H(s,t,...)}{\partial p^i} - \frac{D}{Dt}\frac{\partial H(s,t,...)}{\partial q^i}$$

--- (4.19)

For simplicity, we do not explicitly show the dependence of $g_i$ on $x, x_s, x_t, x_{ss}, x_{tt}, u$ in (4.19). The symbol $\psi_i^{n^*}$ has the meaning $\psi_i^{n^*} \equiv \frac{\partial\psi_i}{\partial n^*} = n_2\frac{\partial\psi_i}{\partial s} + n_1\frac{\partial\psi_i}{\partial t}$. Since $D_j$ is assumed to be non-degenerate, the arc $\partial D_j \cap (\gamma)$ is an oblique part of ($\gamma$), and the values of $\psi_i^{n^*}$ are given by (4.2), and we have

$$\psi_i^{n^*}(s,t) \equiv n_2(s,t)\frac{\partial\psi_i(s,t)}{\partial s} + n_1(s,t)\frac{\partial\psi_i(s,t)}{\partial t} = -2F_{1,i}(s,t,...)$$

--- (4.20)

Consequently, (4.18) is an integral equation of the form

$$\psi_i(s,t) = \psi_{0,i}(s,t) + \iint\limits_{E(s,t)} g_i(\sigma,\tau,\psi(\sigma,\tau),\frac{\partial\psi(\sigma,\tau)}{\partial\sigma},\frac{\partial\psi(\sigma,\tau)}{\partial\tau})\,dA(\sigma,\tau)$$

--- (4.21)

and existence and uniqueness of solutions of integral equations of this type can be shown by proving convergence of a sequence of Picard iterations; this proof is given in appendix B of the present paper. Thus the Hamiltonians are uniquely solvable in the interior of each domain in zone $Z_0$ (cf. definition 4.3). Now, the values $\psi^+_{i(S)}(s,t)$ are known for $(s,t) \in \mathrm{relint}(V_{j+1,j+1}V_{j,j+1})$ from the solution $\psi$ in $\mathrm{int}\,D_j$, and the values $\psi^+_{i(T)}(s,t)$ are known for $(s,t) \in \mathrm{relint}(V_{j+1,j+1}V_{j+1,j+2})$ from the solution $\psi$ in $\mathrm{int}\,D_{j+1}$. The jump conditions of theorem 4.3 then allow the determination of $\psi^-_{i(S)}(s,t)$ for $(s,t) \in \mathrm{relint}(V_{j+1,j+1}V_{j,j+1})$ and $\psi^-_{i(T)}(s,t)$ for $(s,t) \in \mathrm{relint}(V_{j+1,j+1}V_{j+1,j+2})$. Consequently, the Hamiltonian equations for $\psi$ are uniquely solvable as ordinary Goursat-Darboux problems in $\mathrm{int}\,Q_{j,j+1}$ for all j, i.e. $\psi$ is now known in the interior of each domain in zone $Z_1$. Inductively, if $\psi$ is known in the interior of each domain in zone $Z_p$, then the values $\psi^+_{i(S)}(s,t)$ are known for $(s,t) \in \mathrm{relint}(V_{j+1,k+1}V_{j,k+1})$ from the solution $\psi$ in $\mathrm{int}\,Q_{j,k}$, $k - j = p$, and the values $\psi^+_{i(T)}(s,t)$ are known for $(s,t) \in \mathrm{relint}(V_{j+1,k+1}V_{j+1,k+2})$ from the solution $\psi$ in $\mathrm{int}\,Q_{j+1,k+1}$ (trivially, when $Q_{j,k}$ is in zone $Z_p$, $Q_{j+1,k+1}$ is also in zone $Z_p$); the jump conditions then yield the values of $\psi^-_{i(S)}(s,t)$ for $(s,t) \in \mathrm{relint}(V_{j+1,k+1}V_{j,k+1})$ and $\psi^-_{i(T)}(s,t)$ for $(s,t) \in \mathrm{relint}(V_{j+1,k+1}V_{j+1,k+2})$, and then $\psi$ is determined as the unique solution of an ordinary Goursat-Darboux problem in $\mathrm{int}\,Q_{j,k+1}$, thus in the interior of every domain in zone $Z_{p+1}$. The induction is complete and the theorem is proved. ///

<u>Remark 4.1.</u> The treatment presented in sections 3 and 4 can be extended to the more general case in which G is a domain, on the st - plane, containing the origin O in its interior and bounded by a continuous curve $(\gamma)$ with the properties stated below.
Let $G_i$, $i = 1,2,3,4$ denote the parts of G that lie in the corresponding quadrants of the st - plane, and let $(\gamma_i)$, $i = 1,2,3,4$ be the parts of $(\gamma)$ in the same quadrants. (The quadrants of the st - plane are labelled in the standard way, namely the first quadrant is $\{s \ge 0, t \ge 0\}$, the second quadrant is $\{s \le 0, t \ge 0\}$, the third quadrant is $\{s \le 0, t \le 0\}$, and the fourth quadrant is $\{s \geq 0, t \leq 0\}$.) We define the transformed domains $\tilde{G}_i$ and transformed curves $(\tilde{\gamma}_i)$, for i=1, 2, 3, 4, by

$$\tilde{G}_1 := G_1,\ (\tilde{\gamma}_1) := (\gamma_1);\ \tilde{G}_2 = \{(-s,t) : (s,t) \in G_2\},\ (\tilde{\gamma}_2) := \{(-s,t) : (s,t) \in (\gamma_2)\};$$
$$\tilde{G}_3 = \{(-s,-t) : (s,t) \in G_3\},\ (\tilde{\gamma}_3) := \{(-s,-t) : (s,t) \in (\gamma_3)\};$$
$$\tilde{G}_4 = \{(s,-t) : (s,t) \in G_4\},\ (\tilde{\gamma}_4) := \{(s,-t) : (s,t) \in (\gamma_4)\} \qquad \text{--- (4.22)}$$

We postulate that each pair $(G_i,(\gamma_i))$, $i = 1,2,3,4$ should have the same properties as those of the pair (G, (γ)) that were stated in section 2 of this paper.

We consider the problem with state equations

$$\frac{\partial^2 x(s,t)}{\partial s \partial t} = f(s,t,x(s,t),x_s(s,t),x_t(s,t),u(s,t)),\ \text{for } (s,t) \in G;$$
$$x(s,0) = x_1(s) \text{ for } (s,0) \in G;\ x(0,t) = x_2(t) \text{ for } (0,t) \in G \qquad \text{--- (4.23)}$$

and cost functional J, to be minimized, given by

$$J := \iint_G \Phi(s,t,x(s,t),x_s(s,t),x_t(s,t),u(s,t))\,dA(s,t) +$$
$$+ \int_{(\gamma)} \Phi_1(s,t,x(s,t),x_\mu(s,t))\,d\mu(s,t) + \sum_{P \in \mathbf{P}} \Phi_0(P,x(P)) \qquad \text{--- (4.24)}$$

where now the set **P** can contains at most a finite number of points in each part $(\gamma_i), i = 1,2,3,4$ of (γ). The set **P** can be written as a disjoint union $\mathbf{P} = \bigcup_{1 \leq i \leq 4} \mathbf{P}_i$. The functional J can be expressed as a sum $J = \sum_{i=1}^{4} J_i$ with each $J_i$ given by

$$J_i := \iint_{G_i} \Phi(s,t,x(s,t),x_s(s,t),x_t(s,t),u(s,t))\,dA(s,t) +$$

$$+ \int_{(\gamma_i)} \Phi_1(s,t,x(s,t),x_\mu(s,t))\,d\mu(s,t) + \sum_{P \in \mathbf{P}_i} \Phi_0(P,x(P))$$

--- (4.25)

The variation δJ, under a variation δu of the control, can analogously computed as

$$\delta J = \sum_{i=1}^{4} \delta J_i$$

--- (4.26)

and the variation of each $J_i$ can be evaluated by the same techniques as for the standard case of section 2. The co-state ψ can be evaluated in each sub-domain $G_i$ by the same methods as for the standard problem. ///

5. An extremum principle.

Maximum principles of Pontryagin's type have been proved in [BDMO, E1, E2, E3, PS, S, VST] for controlled Goursat-Darboux systems over rectangular domains.
A maximum principle for an optimal control problem has two main components: the Hamiltonian part, i.e. a proof that the co-state satisfies a set of Hamiltonian equations, and the extremum part, i.e. a proof that an optimal control must minimize the Hamiltonian that corresponds to an optimal trajectory. We have already proved the appropriate Hamiltonian equations for the problems of the present paper.
We shall now formulate and prove the extremum part of a maximum principle, akin to Pontryagin's maximum principle. The idea of the variational argument below is based on the approach of Gabasov and Kirillova [GK].

The principal effect of the non-rectangular domain is manifested in the side conditions for the Hamiltonian equations; the proof of the extremum property does not substantially deviate from the proof for rectangular domains, and it is included here primarily for the sake of completeness.

We have:

Theorem 5.1. Assuming the existence of a piecewise continuous optimal control function $u^*(.)$ with corresponding state trajectory $x^*(.,.)$ and co-state trajectory $\psi^*(.,.)$, and assuming that H is uniformly differentiable with respect to the variable u with continuous bounded derivatives $\dfrac{\partial H(s,t,x^*(s,t),x_s^*(s,t),x_t^*(s,t),\psi^*(s,t),u^*(s,t))}{\partial u^k}$ for $(s,t)\in G_R$, we have, for every point (s, t) of $G_R$ that is also a point of continuity of $u^*$,

$$H(s,t,x^*(s,t),x_s^*(s,t),x_t^*(s,t),\psi^*(s,t),u^*(s,t))\le$$
$$\le H(s,t,x^*(s,t),x_s^*(s,t),x_t^*(s,t),\psi^*(s,t),u_1)\ \forall u_1\in \mathbf{U}$$
--- (5.1)

Proof: For every admissible variation $\delta u(.,.)$ of the control function, we denote the finite increment of J by $\Delta J$, thus

$$\Delta J(u^*(.,.),\delta u(.,.)) := J(u^*(.,.)+\delta u(.,.)) - J(u^*(.,.))$$
--- (5.2)

The increment of the Hamiltonian, with respect to a variation in the real variable u(s,t) alone, is denoted by $\Delta_u H$, thus

$$\begin{aligned}&\Delta_u H(s,t,x^*(s,t),x_s^*(s,t),x_t^*(s,t),\psi^*(s,t),u^*(s,t);\delta u(s,t)) :=\\ &= H(s,t,x^*(s,t),x_s^*(s,t),x_t^*(s,t),\psi^*(s,t),u^*(s,t)+\delta u(s,t)) -\\ &- H(s,t,x^*(s,t),x_s^*(s,t),x_t^*(s,t),\psi^*(s,t),u^*(s,t))\end{aligned}$$

--- (5.3)

Because of the definition of J, we have

$$\Delta J(u^*(.,.),\delta u(.,.)) = \delta J(u^*,\delta u) + \mathbf{o}(\| \delta u \|_{L^1(G)})$$

--- (5.4)

where $\delta J(u^*,\delta u)$ signifies $\delta J$ as defined in section 3, evaluated for control function $u^*(.,.)$ and admissible variation $\delta u(.,.)$ . Also, by the differentiability assumption on H, we have

$$\begin{aligned}&\Delta_u H(s,t,x^*(s,t),x_s^*(s,t),x_t^*(s,t),\psi^*(s,t),u^*(s,t);\delta u(s,t)) =\\ &= \frac{\partial H(s,t,x^*(s,t),x_s^*(s,t),x_t^*(s,t),\psi^*(s,t),u^*(s,t))}{\partial u^k}\delta u^k + \mathbf{o}(\| \delta u(s,t) \|)\end{aligned}$$

--- (5.5)

According to the results we have shown in section 3, we have

$$\delta J(u^*,\delta u) = \iint_G \frac{\partial H(s,t,x^*(s,t),x_s^*(s,t),x_t^*(s,t),\psi^*(s,t),u^*(s,t))}{\partial u^k}\delta u^k \, dA(s,t)$$

--- (5.6)

Consequently,

$$\Delta J(u^*(.,.),\delta u(.,.)) =$$
$$= \iint_G \Delta_u H(s,t,x^*(s,t),x_s^*(s,t),x_t^*(s,t),\psi^*(s,t),u^*(s,t);\delta u(s,t))\,dA(s,t) +$$
$$+ \mathbf{o}(\|\delta u\|_{L^1(G)})$$

--- (5.7)

We denote by $B(s,t;\varepsilon)$ the open disk in $IR^2$ with center at (s, t) and radius ε, i.e.

$$B(s,t;\varepsilon) := \{(\sigma,\tau) \in IR^2 : (\sigma - s)^2 + (\tau - t)^2 < \varepsilon^2\}$$

--- (5.8)

We take ε>0 so small that $B(s,t;2\varepsilon) \subseteq G_R$ and $u^*$ is continuous in $B(s,t;2\varepsilon)$. We take as an admissible variation a function $v_\varepsilon$ with the following properties:

$v_\varepsilon$ is continuous in $B(s,t;2\varepsilon)$ ;
$v_\varepsilon$ vanishes outside $B(s,t;\varepsilon)$ ;
$v_\varepsilon(s,t) = u_1 - u^*(s,t)$ , where $u_1$ is an arbitrary but fixed element of **U**.

Then we have

$$\Delta J(u^*(.,.),v_\varepsilon(.,.)) =$$
$$= \iint_G \Delta_u H(s,t,x^*(s,t),x_s^*(s,t),x_t^*(s,t),\psi^*(s,t),u^*(s,t);v_\varepsilon(s,t))\,dA(s,t) +$$
$$+ \mathbf{o}(\varepsilon^2)$$

--- (5.9)

and also

$$\iint\limits_{G} \Delta_u H(s,t,x^*(s,t),x_s^*(s,t),x_t^*(s,t),\psi^*(s,t),u^*(s,t);v_\varepsilon(s,t))\,dA(s,t) =$$

$$= \iint\limits_{B(s,t;\varepsilon)} \Delta_u H(s,t,x^*(s,t),x_s^*(s,t),x_t^*(s,t),\psi^*(s,t),u^*(s,t);v_\varepsilon(s,t))\,dA(s,t) =$$

$$= \pi\varepsilon^2 \Delta_u H(s,t,x^*(s,t),x_s^*(s,t),x_t^*(s,t),\psi^*(s,t),u^*(s,t);u_1 - u^*(s,t)) + \mathbf{o}(\varepsilon^2) \quad \text{--- (5.10)}$$

Therefore

$$\Delta J(u^*(.,.),v_\varepsilon(.,.)) =$$

$$= \pi\varepsilon^2 \Delta_u H(s,t,x^*(s,t),x_s^*(s,t),x_t^*(s,t),\psi^*(s,t),u^*(s,t);u_1 - u^*(s,t)) + \mathbf{o}(\varepsilon^2) \quad \text{--- (5.11)}$$

We use the optimality condition $\Delta J(u^*(.,.),v_\varepsilon(.,.)) \geq 0$ , we divide (5.11) by $\pi\varepsilon^2$ , and we then we let $\varepsilon \to 0^+$, and we obtain

$$\Delta_u H(s,t,x^*(s,t),x_s^*(s,t),x_t^*(s,t),\psi^*(s,t),u^*(s,t);u_1 - u^*(s,t)) \geq 0 \quad \text{--- (5.12)}$$

It is plain that (5.12) is tantamount to (5.1). ///

<u>Remark 5.1.</u> The same extremum principle holds for the more general problem described in remark 4.1. ///

## 6. Application to an inverse problem in the modelling of tsunamis.

The techniques of sections 3, 4, and 5 can be applied to a variety of physical models that involve linear or nonlinear second-order hyperbolic equations in two variables. After changing to characteristic coordinates, we end up with an equation or system of equations in Goursat-Darboux form. Under certain conditions, such problems can fall into the framework of this paper.

As an example, we present a problem that arises in the mathematical modelling of tsunami waves, induced by seismic motions of the bottom of the sea, and also taking into account the earth's rotation. This is only an example, out of many possible examples, of problems to which the methods of the present paper can be applied.

An inverse problem in the modelling of tsunami waves has been treated by optimal control methods in [VST]. Here, we use a model that incorporates additional features and is partially based on the model used in [VS]. Our method of solution also relies on optimal control methods, but using a different procedure than that of [VST]. Our procedure is based on reducing the problem to a Goursat-Darboux equation over a non-rectangular domain, and applying the results of the present paper.

We use the following model:

$$\frac{\partial v(r,t)}{\partial t} - 2\omega w(r,t) + g\frac{\partial y(r,t)}{\partial r} = -\frac{\partial \gamma(r,t)}{\partial t} \; ;$$
$$\frac{\partial w(r,t)}{\partial t} + 2\omega w(r,t) = 0 \; ;$$
$$\frac{\partial y(r,t)}{\partial t} + c\frac{\partial}{\partial r}\left( h(r)\frac{\partial v(r,t)}{\partial r} \right) = 0$$

--- (6.1)

Here, ω is the angular velocity of the basin around a vertical axis z; the equations are written in cylindrical coordinates $(r, \vartheta, z)$ under the assumption of radial symmetry; ζ=η(r) is the equation of the bottom of the sea; v(r, t), w(r, t) are, respectively, the radial and tangential components of the mass velocity of fluid particles; y(r, t) is the elevation of the surface of the water from its relative equilibrium position; z=γ(r, t) describes the vertical motion of the bottom of the sea; c is a proportionality constant; g is the acceleration of gravity.

By taking $\frac{\partial}{\partial t}$ of the first equation in (6.1), we obtain

$$\frac{\partial^2 v(r,t)}{\partial t^2} - 2\omega \frac{\partial w(r,t)}{\partial t} + g \frac{\partial^2 y(r,t)}{\partial r \partial t} = 0$$

--- (6.2)

In view of the second equation in (6.1), we can rewrite (6.2) as

$$\frac{\partial^2 v(r,t)}{\partial t^2} + 4\omega^2 v(r,t) + g \frac{\partial^2 y(r,t)}{\partial r \partial t} = 0$$

--- (6.3)

Differentiation of the third equation in (6.1) with respect to r gives

$$\frac{\partial^2 y(r,t)}{\partial r \partial t} + c \frac{\partial^2}{\partial r^2}\big(h(r)v(s,r)\big) = -\frac{\partial^2 \gamma(r,t)}{\partial r \partial t}$$

--- (6.4)

Elimination of the term $\frac{\partial^2 y(r,t)}{\partial r \partial t}$ between (6.3) and (6.4) yields

$$\frac{\partial^2 v(r,t)}{\partial t^2} + 4\omega^2 v(r,t) - gc \frac{\partial^2}{\partial r^2}\big(h(r)v(r,t)\big) = g \frac{\partial^2 \gamma(r,t)}{\partial r \partial t}$$

--- (6.5)

We set

$$u(r,t) := \frac{\partial^2 \gamma(r,t)}{\partial r \partial t}$$

--- (6.6)

We denote by $s, t'$ the characteristic coordinates for (6.5). These are found as

$$s = s_0 + t + \frac{1}{\sqrt{gc}}\beta(r) \ ; \ \ t' = t_0' + t - \frac{1}{\sqrt{gc}}\beta(r) \ ; \ \ \beta(r) := \int_0^r \frac{1}{\sqrt{h(\rho)}}\, d\rho$$

--- (6.7)

(Naturally, we assume that the origin of coordinates is chosen so that $h(r) > 0$ for all values of r.)

The observation domain G is taken as

$$G := \{(r, t): r_1 \le r \le r_2 \, , \, t_1 \le t \le t_2\} \qquad \text{--- (6.8)}$$

where $r_1, r_2, t_1, t_2$ satisfy the condition

$$\beta(r_2) - \beta(r_1) = \sqrt{gc}\,(t_2 - t_1) \qquad \text{--- (6.9)}$$

Then, for a suitable choice of $s_0, t_0'$, the domain G is described, in characteristic coordinates, by

$$G := \{(s, t'): -A \le s + t' \le A \, , \, -A \le s - t' \le A\} \qquad \text{--- (6.10)}$$

The values of A, $s_0, t_0'$ are determined as solutions of the system

$$\begin{aligned} A &= s_0 + t_0' + 2t_2 = s_0 - t_0' + \frac{2}{\sqrt{gc}}\beta(r_2) = \\ &= -(s_0 + t_0') - 2t_1 = -(s_0 - t_0') - \frac{2}{\sqrt{gc}}\beta(r_1) \end{aligned} \qquad \text{--- (6.11)}$$

Clearly, the first 3 out of the 4 equations in (6.11) suffice for the determination of A, $s_0, t_0'$; the fourth equation is consistent with the other 3 because of the condition (6.9).

The partial differential equation for v becomes

$$\frac{\partial^2 v(s,t')}{\partial s \partial t'} = \frac{1}{2}\left(\frac{\partial v(s,t')}{\partial t'} - \frac{\partial v(s,t')}{\partial s}\right)\frac{\sqrt{gc}\,h'(r)}{\sqrt{h(r)}} + \left(\frac{gch''(r)}{4} + \omega^2\right)v(s,t') + \frac{gc}{4}u(s,t') \quad \text{--- (6.12)}$$

In (6.12), r is expressed as a function of the characteristic coordinates as

$$r \equiv r(s,t') = \beta^{-1}\left(\frac{\sqrt{gc}}{2}(s - t' - s_0 + t_0')\right) \quad \text{--- (6.13)}$$

where $\beta^{-1}$ is the inverse function of β, and it is well defined since $\frac{d\beta(r)}{dr} = \frac{1}{\sqrt{h(r)}} > 0$.

We examine the problem of estimating u from observations $\tilde{v}(s,t')$ of the function v over the domain G. We formulate this problem as a problem of optimal control, and we seek to minimize a functional J given by

$$J = \iint_G F(s,t',|\,v(s,t') - \tilde{v}(s,t')\,|^2, u(s,t'))\,ds\,dt' \quad \text{--- (6.14)}$$

For example, we may take

$$F(s,t',|\,v(s,t') - \tilde{v}(s,t')\,|^2, u(s,t')) = |\,v(s,t') - \tilde{v}(s,t')\,|^2 + \lambda\,|\,u(s,t')\,|^2 \quad \text{--- (6.15)}$$

where λ is a regularization parameter.

We note that this problem is different from the related problem in [VST]. Here, we use only observations of v, but not observations of y, while y was required in [VST]; also, here we estimate $\frac{\partial^2 \gamma(r,t)}{\partial r \partial t}$, whereas the estimated quantity in [VST] was $\frac{\partial \gamma}{\partial t}$.

The domain G is made up of 4 parts $G_i$, $1 \le i \le 4$, where each $G_i$ is the part of G that lies in the i-th quadrant of the $st'$ - plane. Each of those 4 parts satisfies the conditions used in the above sections of the present paper, and necessary optimality conditions have the form of the Hamiltonian equation , the extremum principle of section 5, and the side conditions of section 4 on each part $(\gamma_i)$ of the boundary of G in each of the 4 quadrants.

The numerical realization of the optimality conditions of sections 3-5, for this specific problem of tsunami modelling, is beyond the scope of the present paper.

Appendix A: integration and differentiation formulae.

In this section, we present a number of results that are needed for the manipulation of integrals and derivatives arising in the derivation of the Hamiltonian equations.

We shall use the notation and terminology of section 2 of this paper. In addition, we define

$$E(s,t) := \{(\sigma,\tau) \in G : (s,t) \in W(\sigma,\tau)\} \qquad \text{--- (A.1)}$$

Consequently, for functions $\varphi(s,t,\sigma,\tau)$ that are integrable over $G \times G$, we have, by Fubini's theorem,

$$\iint_{G} \iint_{W(s,t)} \varphi(s,t,\sigma,\tau)\, dA(\sigma,\tau)\, dA(s,t) = \iint_{G} \iint_{E(s,t)} \varphi(\sigma,\tau,s,t)\, dA(\sigma,\tau)\, dA(s,t) \qquad \text{--- (A.2)}$$

For every point $(s,t) \in \operatorname{int} G$, we denote by $A_t$, $B_s$ the points of intersection of $(\gamma)$ with the straight lines through (s, t) parallel to the s-axis and to the t-axis, respectively; thus

$$A_t = (\vartheta_1(t), t),\ B_s = (s, \vartheta_2(s)) \qquad \text{--- (A.3)}$$

Therefore,

$$(\sigma,\tau) \in \gamma(s,t) \Leftrightarrow (s,t) \in W(\sigma,\tau) \qquad \text{--- (A.4)}$$

and Fubini's theorem for functions $\varphi(s,t,\sigma,\tau)$ that are integrable over $(\gamma) \times G$ takes the form

$$\int_{(\gamma)} \iint_{W(s,t)} \varphi(s,t,\sigma,\tau)\, dA(\sigma,\tau)\, d\mu(s,t) = \iint_{G} \int_{\gamma(s,t)} \varphi(\sigma,\tau,s,t)\, d\mu(\sigma,\tau)\, dA(s,t) \qquad \text{--- (A.5)}$$

For the sake of notational consistency, the point (s, t) will also be denoted by $P_{st}$ . When the function $\varphi$ is absolutely continuous (in the sense of absolute continuity for functions of two variables), we have, for almost all (s, t) in G,

$$\frac{\partial}{\partial s}\iint_{E(s,t)} \varphi(s,t,\sigma,\tau)\,dA(\sigma,\tau) = -\int_{P_{st}}^{B_s} \varphi(s,t,s,\tau)\,d\tau + \iint_{E(s,t)} \varphi_s(s,t,\sigma,\tau)\,dA(\sigma,\tau)$$

--- (A.6)

$$\frac{\partial}{\partial t}\iint_{E(s,t)} \varphi(s,t,\sigma,\tau)\,dA(\sigma,\tau) = -\int_{P_{st}}^{A_t} \varphi(s,t,\sigma,t)\,d\sigma + \iint_{E(s,t)} \varphi_t(s,t,\sigma,\tau)\,dA(\sigma,\tau)$$

--- (A.7)

$$\frac{\partial^2}{\partial s\partial t}\iint_{E(s,t)} \varphi(s,t,\sigma,\tau)\,dA(\sigma,\tau) = \varphi(s,t,s,t) - -\int_{P_{st}}^{A_t} \varphi_s(s,t,\sigma,t)\,d\sigma - \int_{P_{st}}^{B_s} \varphi_t(s,t,s,\tau)\,d\tau +$$

$$+ \iint_{E(s,t)} \varphi_{st}(s,t,\sigma,\tau)\,dA(\sigma,\tau)$$

--- (A.8)

$$\frac{\partial}{\partial s}\int_{\gamma(s,t)} \varphi(s,t,\sigma,\tau)\,d\mu(\sigma,\tau) = -\frac{\varphi(s,t,B_s)}{n_2(B_s)} + \int_{\gamma(s,t)} \frac{\partial\varphi(s,t,\sigma,\tau)}{\partial s}\,d\mu(\sigma,\tau)$$

--- (A.9)

$$\frac{\partial}{\partial t}\int_{\gamma(s,t)} \varphi(s,t,\sigma,\tau)\,d\mu(\sigma,\tau) = -\frac{\varphi(s,t,A_t)}{n_1(A_t)} + \int_{\gamma(s,t)} \frac{\partial\varphi(s,t,\sigma,\tau)}{\partial t}\,d\mu(\sigma,\tau)$$

--- (A.10)

$$\frac{\partial^2}{\partial s \partial t} \int_{\gamma(s,t)} \varphi(s,t,\sigma,\tau)\, d\mu(\sigma,\tau) = -\frac{\varphi_s(s,t,A_t)}{n_1(A_t)} - \frac{\varphi_t(s,t,B_s)}{n_2(B_s)} +$$
$$+ \int_{\gamma(s,t)} \frac{\partial^2 \varphi(s,t,\sigma,\tau)}{\partial s \partial t} d\mu(\sigma,\tau)$$

--- (A.11)

We note that, for $(s,t) \in G_R$, the quantities $n_1(A_t), n_2(B_s)$ are both nonzero, because of the geometric meaning of the points $A_t$ and $B_s$.

Next, we prove one formula out of the group (A.6), (A.7), and (A.8). As a representative case, we prove (A.6), and the proofs of the remaining formulae will be similar. For simplicity of the exposition, and without loss of generality, we present the proof for one of the possible cases, namely for the case in which the arc $A_t B_s \subseteq (\gamma)$ has the representation $\tau = \vartheta_2(\sigma)$. We have

$$\iint_{E(s,t)} \varphi(s,t,\sigma,\tau)\, dA(\sigma,\tau) = \int_{\sigma=s}^{s(A_t)} \int_{\tau=t}^{\vartheta_2(\sigma)} \varphi(s,t,\sigma,\tau)\, d\tau\, d\sigma$$

--- (A.12)

and consequently

$$\frac{\partial}{\partial s} \iint_{E(s,t)} \varphi(s,t,\sigma,\tau)\, dA(\sigma,\tau) = -\int_{\tau=t}^{\vartheta_2(s)} \varphi(s,t,s,\tau)\, d\tau + \int_{\sigma=s}^{s(A_t)} \int_{\tau=t}^{\vartheta_2(\sigma)} \varphi_s(s,t,\sigma,\tau)\, d\tau\, d\sigma =$$
$$= -\int_{P_{st}}^{B_s} \varphi(s,t,s,\tau)\, d\tau + \iint_{E(s,t)} \varphi_s(s,t,\sigma,\tau)\, dA(\sigma,\tau)$$

--- (A.13)

Similarly we prove formula (A.9) out of the group of formulae (A.9), (A.10), (A.11), the proofs of the remaining formulae being similar.
When the arc $\gamma(s,t)$ can be represented both in the form $\sigma = \vartheta_1(\tau)$ and in the form $\tau = \vartheta_2(\sigma)$, the proof becomes somewhat simplified, and we present that simpler proof first. We have

$$\int_{\gamma(s,t)} \varphi(s,t,\sigma,\tau)\,d\mu(\sigma,\tau) = \int_{\tau=t}^{\vartheta_2(s)} \varphi(s,t,\vartheta_1(\tau),\tau)\sqrt{1+(\vartheta_1'(\tau))^2}\,d\tau$$

--- (A.14)

and consequently

$$\frac{\partial}{\partial s}\int_{\gamma(s,t)} \varphi(s,t,\sigma,\tau)\,d\mu(\sigma,\tau) = \vartheta_2'(s)\varphi(s,t,s,\vartheta_2(s))\sqrt{1+(\vartheta_1'(\tau))^2}\,\Bigg|_{\tau=\vartheta_2(s)} +$$

$$+ \int_{\gamma(s,t)} \varphi_s(s,t,\sigma,\tau)\,d\mu(\sigma,\tau) =$$

$$= \vartheta_2'(s)\varphi(s,t,s,\vartheta_2(s))\sqrt{1+\left(\frac{1}{\vartheta_2'(s)}\right)^2} + \int_{\gamma(s,t)} \varphi_s(s,t,\sigma,\tau)\,d\mu(\sigma,\tau) =$$

$$-\sqrt{1+(\vartheta_2'(s))^2}\,\varphi(s,t,s,\vartheta_2(s)) + \int_{\gamma(s,t)} \varphi_s(s,t,\sigma,\tau)\,d\mu(\sigma,\tau) =$$

$$= -\frac{\varphi(s,t,B_s)}{n_2(B_s)} + \int_{\gamma(s,t)} \varphi_s(s,t,\sigma,\tau)\,d\mu(\sigma,\tau)$$

--- (A.15)

In the general case, we denote by $\mu_t$, $\mu_s$ the values of the parameter $\mu$ that correspond to the points $A_t$, $B_s$, respectively. Then

$$\int_{\gamma(s,t)} \varphi(s,t,\sigma,\tau)\,d\mu(\sigma,\tau) = \int_{\mu_t}^{\mu_s} \varphi(s,t,L_1(\mu),L_2(\mu))\,d\mu$$

--- (A.16)

and, consequently,

$$\frac{\partial}{\partial s}\int_{\gamma(s,t)} \varphi(s,t,\sigma,\tau)\,d\mu(\sigma,\tau) = \varphi(s,t,L_1(\mu_s),L_2(\mu_s))\frac{\partial \mu_s}{\partial s} + \int_{\mu_t}^{\mu_s} \frac{\partial \varphi(s,t,L_1(\mu),L_2(\mu))}{\partial s}\,d\mu$$

--- (A.17)

It follows from the geometric definition of $A_t$ and $B_s$ that, near $B_s$, the curve (γ) has the representation $\tau = \vartheta_2(\sigma)$. When t is kept constant and s changes to $s + \delta s$, the corresponding point $(s, \tau)$ on $(\gamma)$, with $\tau = \vartheta_2(s)$, changes to $(s + \delta s, \tau + \delta\tau)$ where $\delta\tau = \vartheta_2{}'(s)\delta s + o(\delta s)$, and $\mu_s$ changes to $\mu_s + \delta\mu$ where $\delta\mu = -\sqrt{1 + (\vartheta_2{}'(s))^2}\,\delta s + o(\delta s)$. Therefore,

$$\frac{\partial \mu_s}{\partial s} = -\sqrt{1 + (\vartheta_2{}'(s))^2} = -\frac{1}{n_2(B_s)} \qquad \text{--- (A.18)}$$

and then (A.9) follows from (A.17) and (A.18).

Appendix B: solvability of certain integral equations.

We examine the system of integral equations of the form (4.21), which we rewrite here for the reader's convenience:

$$\psi_i(s,t) = \psi_{0,i}(s,t) + \iint\limits_{E(s,t)} g_i(\sigma,\tau,\psi(\sigma,\tau),\frac{\partial\psi(\sigma,\tau)}{\partial\tau},\frac{\partial\psi(\sigma,\tau)}{\partial\sigma})\,dA(\sigma,\tau) \quad \text{--- (B.1)}$$

The conditions on E(s, t) and the domain $D_j$ over which a solution of (B.1) is sought are as stated in section 4.

We note that the method of proof of this section can also be used for a more general integral equation, in which the functions $g_i$ inside the integral also depend on s and t, i.e. the case $g_i = g_i(s,t,\sigma,\tau,\psi(\sigma,\tau),\frac{\partial\psi(\sigma,\tau)}{\partial\tau},\frac{\partial\psi(\sigma,\tau)}{\partial\sigma})$. However, this extension is not needed for the present paper.

The simplest existence and uniqueness result is achieved through the standard method of Picard iterations.

We assume that each $g_i$ is continuous with respect to all its arguments, and satisfies a Lipschitz condition

$$|g_i(\sigma,\tau,\Psi_1,P_1,Q_1) - g_i(\sigma,\tau,\Psi_2,P_2,Q_2)| \le$$
$$\le L\left[|\Psi_1 - \Psi_2|_{2,IR^n} + |P_1 - P_2|_{2,IR^n} + |Q_1 - Q_2|_{2,IR^n}\right] \quad \text{--- (B.2)}$$

The symbol $|\cdot|_{2,IR^n}$ denotes the Euclidean norm on $IR^n$. Since no assumption is made on the magnitude of L, it is clear that any other norm can be used in the Lipschitz conditions.

We seek a solution of (B.1) in the space $C_1(D_j;IR^n)$ of $IR^n$ - valued functions that are continuous in $D_j$ and have first partial derivatives in $\operatorname{int} D_j$ which are continuous in $D_j$.

By differentiating (B.1) with respect to s and t, and using the differentiation formulae of appendix A, we obtain

$$\frac{\partial\psi_i(s,t)}{\partial s} = \frac{\partial\psi_{0,i}(s,t)}{\partial s} - \int_{P_{st}}^{B_s} g_i(s,\tau,\psi(s,\tau),\frac{\partial\psi(s,\tau)}{\partial s},\frac{\partial\psi(s,\tau)}{\partial\tau})\,d\tau\ ;$$

$$\frac{\partial\psi_i(s,t)}{\partial t} = \frac{\partial\psi_{0,i}(s,t)}{\partial t} - \int_{P_{st}}^{A_t} g_i(\sigma,t,\psi(\sigma,t),\frac{\partial\psi(\sigma,t)}{\partial\sigma},\frac{\partial\psi(\sigma,t)}{\partial t})\,d\sigma$$

--- (B.3)

We set

$$P_{(i)}(s,t) := \frac{\partial\psi_i(s,t)}{\partial s}\ ,\ Q_{(i)}(s,t) := \frac{\partial\psi_i(s,t)}{\partial t}$$

--- (B.4)

We replace (B.1) by the system

$$\psi_i(s,t) = \psi_{0,i}(s,t) + \iint\limits_{E(s,t)} g_i(\sigma,\tau,\psi(s,t),P(s,t),Q(s,t))\,dA(s,t)\ ;$$

$$P_{(i)}(s,t) = \frac{\partial\psi_{0,i}(s,t)}{\partial s} - \int_{P_{st}}^{B_s} g_i(s,\tau,\psi(s,\tau),P(s,\tau),Q(s,\tau))\,d\tau\ ;$$

$$Q_{(i)}(s,t) = \frac{\partial\psi_{0,i}(s,t)}{\partial s} - \int_{P_{st}}^{A_t} g_i(\sigma,t,\psi(\sigma,t),P(\sigma,t),Q(\sigma,t))\,d\sigma$$

--- (B.5)

It can be shown, by using the differentiation formulae (B.2), that the problem (B.5) is equivalent to the original problem (B.1).

For each $\rho \geq 0$ , we define the norms

$$\|\psi_i\|_\rho := \max_{(s,t)\in D_j} \exp(-\rho(a-s+b-t))\,|\psi_i(s,t)|\ ;$$

$$\|\psi\|_\rho := \max_{(s,t)\in D_j} \exp(-\rho(a-s+b-t))\,|\psi(s,t)|_{2,IR^n}\ ;$$

$$\|(\psi_i,P_{(i)},Q_{(i)})\|_\rho := \|\psi_i\|_\rho + \|P_{(i)}\|_\rho + \|Q_{(i)}\|_\rho\ ;$$

$$\|(\psi,P,Q)\|_\rho := \|\psi\|_\rho + \|P\|_\rho + \|Q\|_\rho$$

--- (B.6)

(The norms $\| P_{(i)} \|_\rho$ , $\| Q_{(i)} \|_\rho$ , $\| P \|_\rho$ , $\| Q \|_\rho$ are defined in the same way as the norms $\| \psi_i \|_\rho$ , $\| \psi \|_\rho$ .)

We define the operators **Ψ**, **P**, **Q**, **S** by

$$\left(\boldsymbol{\Psi}_i(\psi,P,Q)\right)(s,t) := \psi_{0,i}(s,t) + \iint\limits_{E(s,t)} g_i(\sigma,\tau,\psi(s,t),P(s,t),Q(s,t))\,dA(s,t)\ ;$$
$$\left(\mathbf{P}_i(\psi,P,Q)\right)(s,t) := \frac{\partial\psi_{0,i}(s,t)}{\partial s} - \int_{P_{st}}^{B_s} g_i(s,\tau,\psi(s,\tau),P(s,\tau),Q(s,\tau))\,d\tau\ ;$$
$$\left(\mathbf{Q}_i(\psi,P,Q)\right)(s,t) := \frac{\partial\psi_{0,i}(s,t)}{\partial s} - \int_{P_{st}}^{A_t} g_i(\sigma,t,\psi(\sigma,t),P(\sigma,t),Q(\sigma,t))\,d\sigma\ ;$$
$$\left(\mathbf{S}_i(\psi,P,Q)\right)(s,t) := \left(\left(\boldsymbol{\Psi}_i(\psi,P,Q)\right)(s,t),\left(\mathbf{P}_i(\psi,P,Q)\right)(s,t),\left(\mathbf{Q}_i(\psi,P,Q)\right)(s,t)\right);$$
$$\left(\boldsymbol{\Psi}(\psi,P,Q)\right)(s,t) := \left(\left(\boldsymbol{\Psi}_i(\psi,P,Q)\right)(s,t) : 1 \le i \le n\right);$$
$$\left(\mathbf{P}(\psi,P,Q)\right)(s,t) := \left(\left(\mathbf{P}_i(\psi,P,Q)\right)(s,t) : 1 \le i \le n\right);$$
$$\left(\mathbf{Q}(\psi,P,Q)\right)(s,t) := \left(\left(\mathbf{Q}_i(\psi,P,Q)\right)(s,t) : 1 \le i \le n\right);$$
$$\left(\mathbf{S}(\psi,P,Q)\right)(s,t) := \left(\left(\boldsymbol{\Psi}(\psi,P,Q)\right)(s,t),\left(\mathbf{P}(\psi,P,Q)\right)(s,t),\left(\mathbf{Q}(\psi,P,Q)\right)(s,t)\right)$$

--- (B.7)

For two sets of functions, $\left(\psi_{i[1]},P_{(i)[1]},Q_{(i)[1]}\right)$, $\left(\psi_{i[2]},P_{(i)[2]},Q_{(i)[2]}\right)$, we have the estimates

$$\exp(-\rho(a-s+b-t))\,|\,(\mathbf{\Psi}_i(\psi_{[1]},P_{[1]},Q_{[1]}))(s,t)-(\mathbf{\Psi}_i(\psi_{[2]},P_{[2]},Q_{[2]}))(s,t)\,|\leq$$

$$\leq \exp(-\rho(a-s+b-t)) \iint_{E(s,t)} |\,g_i(\sigma,\tau,\psi_{[1]}(\sigma,\tau),P_{[1]}(\sigma,\tau),Q_{[1]}(\sigma,\tau))-$$

$$-g_i(\sigma,\tau,\psi_{[2]}(\sigma,\tau),P_{[2]}(\sigma,\tau),Q_{[2]}(\sigma,\tau))\,|\,dA(\sigma,\tau)\leq$$

$$\leq \exp(-\rho(a-s+b-t)) \iint_{E(s,t)} \exp(\rho(a-\sigma+b-\tau))\,L\,(\|\,\psi_{[1]}-\psi_{[2]}\,\|_\rho +$$

$$+\|\,P_{[1]}-P_{[2]}\,\|_\rho+\|\,Q_{[1]}-Q_{[2]}\,\|_\rho)\,dA(\sigma,\tau)=$$

$$=\exp(-\rho(a-s+b-t))\,L\,(\|\,\psi_{[1]}-\psi_{[2]}\,\|_\rho+\|\,P_{[1]}-P_{[2]}\,\|_\rho+\|\,Q_{[1]}-Q_{[2]}\,\|_\rho)\cdot$$

$$\cdot \iint_{E(s,t)} \exp(\rho(a-\sigma+b-\tau))\,dA(\sigma,\tau)\leq$$

$$\leq \exp(-\rho(a-s+b-t))\,L\,(\|\,\psi_{[1]}-\psi_{[2]}\,\|_\rho+\|\,P_{[1]}-P_{[2]}\,\|_\rho+\|\,Q_{[1]}-Q_{[2]}\,\|_\rho)\cdot$$

$$\cdot\int_s^a \int_t^b \exp(\rho(a-\sigma+b-\tau))\,d\tau\,d\sigma=$$

$$=\frac{(1-\exp(-\rho(s-a)))(1-\exp(\rho(t-b)))}{\rho^2}\cdot$$

$$\cdot L\,(\|\,\psi_{[1]}-\psi_{[2]}\,\|_\rho+\|\,P_{[1]}-P_{[2]}\,\|_\rho+\|\,Q_{[1]}-Q_{[2]}\,\|_\rho)\leq$$

$$\leq\frac{(1-\exp(-\rho a))(1-\exp(-\rho b))}{\rho^2}\,L\,(\|\,\psi_{[1]}-\psi_{[2]}\,\|_\rho+\|\,P_{[1]}-P_{[2]}\,\|_\rho+\|\,Q_{[1]}-Q_{[2]}\,\|_\rho)$$

--- (B.8)

$$\exp(-\rho(a-s+b-t))\,|\,\big(\mathbf{P}_i(\psi_{[1]},P_{[1]},Q_{[1]})\big)(s,t)-\big(\mathbf{P}_i(\psi_{[2]},P_{[2]},Q_{[2]})\big)(s,t)\,|\le$$
$$\le \exp(-\rho(a-s+b-t))\int_{P_{st}}^{B_s} |\,g_i(s,\tau,\psi_{[1]}(s,\tau),P_{[1]}(s,\tau),Q_{[1]}(s,\tau))-$$
$$-g_i(s,\tau,\psi_{[2]}(s,\tau),P_{[2]}(s,\tau),Q_{[2]}(s,\tau))\,|\,d\tau\le$$
$$\le \exp(-\rho(a-s+b-t))\int_{P_{st}}^{B_s} \exp(\rho(b-\tau))\,L\,(\|\,\psi_{[1]}-\psi_{[2]}\,\|_\rho +$$
$$+\|\,P_{[1]}-P_{[2]}\,\|_\rho+\|\,Q_{[1]}-Q_{[2]}\,\|_\rho)\,d\tau=$$
$$=\exp(-\rho(a-s+b-t))\,L\,(\|\,\psi_{[1]}-\psi_{[2]}\,\|_\rho+\|\,P_{[1]}-P_{[2]}\,\|_\rho+\|\,Q_{[1]}-Q_{[2]}\,\|_\rho)\cdot$$
$$\cdot\int_{P_{st}}^{B_s} \exp(\rho(b-\tau))\,d\tau\le$$
$$\le \exp(-\rho(a-s+b-t))\,L\,(\|\,\psi_{[1]}-\psi_{[2]}\,\|_\rho+\|\,P_{[1]}-P_{[2]}\,\|_\rho+\|\,Q_{[1]}-Q_{[2]}\,\|_\rho)\cdot$$
$$\cdot\int_t^b \exp(\rho(b-\tau))\,d\tau\, d\sigma=$$
$$=\frac{(1-\exp(-\rho(s-a)))(1-\exp(\rho(t-b)))}{\rho^2}\cdot$$
$$\cdot\,L\,(\|\,\psi_{[1]}-\psi_{[2]}\,\|_\rho+\|\,P_{[1]}-P_{[2]}\,\|_\rho+\|\,Q_{[1]}-Q_{[2]}\,\|_\rho)\le$$
$$\le\frac{(1-\exp(-\rho a))(1-\exp(-\rho b))}{\rho^2}\,L\,(\|\,\psi_{[1]}-\psi_{[2]}\,\|_\rho+\|\,P_{[1]}-P_{[2]}\,\|_\rho+\|\,Q_{[1]}-Q_{[2]}\,\|_\rho)$$

--- (B.9)

with a similar estimate for the operator **Q**.

Since $\lim_{\rho\to\infty}\frac{(1-\exp(-\rho a))(1-\exp(-\rho b))}{\rho^2}=0$, it follows that, for ρ sufficiently large, the operator **S** is a contraction relative to the norm $\|\,(\psi,P,Q)\,\|_\rho$, and consequently it has a fixed point which is obtained as the uniform limit of the Picard iterations corresponding to the system (B.5). Thus, if we take an arbitrary triple $\big(\psi_{[0]}\ ,\ P_{[0]}\ ,\ Q_{[0]}\big)$ in $\big(C(D_j;IR^n)\big)^3$, and define a sequence of iterations inductively by

$$\big(\psi_{[n+1]},P_{[n+1]},Q_{[n+1]}\big)=\mathbf{S}\big(\psi_{[n]},P_{[n]},Q_{[n]}\big)$$

--- (B.10)

then the following limit exists, in the norm $\|\,(\psi,P,Q)\,\|_\rho$, thus in the topology of uniform convergence in $\big(C(D_j;IR^n)\big)^3$ :

$$\lim_{n\to\infty}\left(\psi_{[n]},P_{[n]},Q_{[n]}\right)\equiv\left(\psi_{[\infty]},P_{[\infty]},Q_{[\infty]}\right)$$

--- (B.11)

and the function $\psi_{[\infty]}$ is the unique solution of (B.1), whereas

$P_{[\infty]}=\dfrac{\partial\psi_{[\infty]}}{\partial s}$ , $Q_{[\infty]}=\dfrac{\partial\psi_{[\infty]}}{\partial t}$ .